\documentclass[a4paper]{amsart}

\usepackage[T1]{fontenc}
\usepackage[latin1]{inputenc}
\usepackage{amscd}
\usepackage{times}
\usepackage{epsfig}
\makeatletter

\theoremstyle{plain}    
\newtheorem{thm}{Theorem}[section]
\newtheorem{defn}[thm]{Definition}

\newtheorem{question}[thm]{Question}
\newtheorem{summary}[thm]{Summary}
\numberwithin{equation}{section} 
\numberwithin{figure}{section} 
\theoremstyle{plain}    
\newtheorem{cor}[thm]{Corollary} 
\newtheorem{lem}[thm]{Lemma} 
\theoremstyle{plain}    
\newtheorem{prop}[thm]{Proposition} 
\newtheorem{fact}[thm]{Fact} 
\theoremstyle{remark}
\newtheorem{claim}[thm]{Claim} 
\newtheorem{rem}[thm]{Remark}
\theoremstyle{remark}    
\newtheorem{notation}[thm]{Notation} 
\newtheorem{assumption}[thm]{Assumption}

\def\factor#1.#2.{\left. \raise 2pt\hbox{$#1$} \right/\hskip -2pt\raise -2pt\hbox{$#2$}}

\usepackage[arrow,matrix,curve,ps]{xy}

\newcommand{\ilabel}[1]{\newcounter{#1}\setcounter{#1}{\value{enumi}}}
\newcommand{\iref}[1]{\setcounter{enumi}{\value{#1}}\labelenumi}

\newcommand{\C}{\mbox{$\mathbb{C}$}}

\renewcommand{\O}{\mbox{$\mathcal{O}$}}
\renewcommand{\P}{\mathbb{P}}

\DeclareMathOperator{\Aut}{Aut}
\DeclareMathOperator{\BlowUp}{BlowUp}

\DeclareMathOperator{\Univrc}{{Univ^{rc}}}

\DeclareMathOperator{\Chow}{Chow}
\DeclareMathOperator{\codim}{codim}

\DeclareMathOperator{\Ext}{Ext}

\DeclareMathOperator{\Hom}{Hom}

\DeclareMathOperator{\Image}{Image}
\DeclareMathOperator{\Jet}{Jet}
\DeclareMathOperator{\locus}{locus}

\DeclareMathOperator{\Pic}{Pic}
\DeclareMathOperator{\Prolong}{Prolong}
\DeclareMathOperator{\rank}{rank}
\DeclareMathOperator{\RatCurves}{RatCurves^n}

\DeclareMathOperator{\red}{red}

\DeclareMathOperator{\Sing}{Sing}

\makeatother

\begin{document}

\title{Lines on complex contact manifolds IIb}

\begin{abstract}
Let $X$ be a complex-projective contact manifold with $b_2(X)=1$. It
has long been conjectured that $X$ should then be
rational-homogeneous, or equivalently, that there exists an embedding
$X \to \P^n$ whose image contains lines. 

We show that $X$ is covered by a compact family of rational curves,
called ``contact lines'' that behave very much like the lines on the
rational homogeneous examples: if $x \in X$ is a general point, then
all contact lines through $x$ are smooth, no two of them share a
common tangent direction at $x$, and the union of all contact lines
through $x$ forms a cone over an irreducible, smooth base. As a
corollary, we obtain that the tangent bundle of $X$ is stable.
\end{abstract}

\date{\today}

\author{Stefan Kebekus}



\thanks{The author was supported by a Heisenberg-Fellowship and by the
  program ``Globale Methoden in der Komplexen Analysis'' of the
  Deutsche Forschungsgemeinschaft.}

\address{Stefan Kebekus, Mathematisches Institut, Universität zu Köln, Weyertal~86--90, 50931~Köln, Germany}
\email{\href{mailto:stefan.kebekus@math.uni-koeln.de}{stefan.kebekus@math.uni-koeln.de}}
\urladdr{\href{http://www.mi.uni-koeln.de/~kebekus}{http://www.mi.uni-koeln.de/$\sim$kebekus}}

\maketitle
\tableofcontents

\section{Introduction}

Motivated by questions coming from Riemannian geometry, complex
contact manifolds have received considerable attention during the last
years.  The link between complex and Riemannian geometry is given by
the twistor space construction: twistor spaces over Riemannian
manifolds with quaternion-Kähler holonomy group are complex contact
manifolds. As twistor spaces are covered by rational curves, much of
the research is centered about the geometry of rational curves on the
contact spaces.

\subsection{Setup and Statement of the main result}

Throughout the present paper, we maintain the assumptions and
notational conventions of the first part \cite{Kebekus01} of this
article. In particular, we refer to \cite{Kebekus01}, and the
references therein, for an introduction to contact manifolds and to
the parameter spaces which we will use freely throughout.

In brief, we assume throughout that $X$ is a complex projective
manifold of dimension $\dim X = 2n+1$ which carries a contact
structure. This structure is given by a vector bundle sequence
\begin{equation}\label{eq:contact-sequence}
  \begin{CD}
    0 @>>> F @>>> T_X @>{\theta}>> L @>>> 0
  \end{CD}
\end{equation}
where $F$ is a subbundle of corank 1 and where the skew-symmetric
O'Neill tensor
$$
N : F\otimes F \to L,
$$
which is associated with the Lie-Bracket, is non-degenerate at
every point of $X$.

Because contact manifolds with $b_2(X)>1$ were completely described in
\cite{KPSW00}, we consider only the case where $b_2(X)=1$. We will
also assume that $X$ is not isomorphic to the projective space
$\P^{2n+1}$. By \cite[Sect.~2.3]{Kebekus01}, these assumptions imply
that we can find a compact irreducible component $H \subset
\RatCurves(X)$ of the space of rational curves on $X$ such that the
intersection of $L$ with the curves associated with $H$ is one. Curves
that are associated with points of $H$ are called ``contact lines''.
For a point $x \in X$, consider the varieties
$$
H_x := \{ \ell \in H \,|\, x \in \ell\}
$$
and
$$
\locus (H_x) := \bigcup_{\ell \in H_x} \ell.
$$
The main result of this paper is the following:

\begin{thm}\label{thm:main}
  Let $X$ be a complex-projective contact manifold with $b_2(X)=1$ and
  assume $X \not \cong \P^{\dim X}$. Let $H \subset \RatCurves(X)$ be
  an irreducible component which parameterizes contact lines. Then
  $\locus(H_x)$ is isomorphic to a projective cone over a smooth,
  irreducible base. Further,
  \begin{enumerate}
  \item \ilabel{l1} all contact lines that contain $x$ are smooth,
    
  \item \ilabel{l2}the space $H_x$ is irreducible,
  
  \item \ilabel{l3} if $\ell_1$ and $\ell_2$ are any two contact lines
    through $x$, then $T_{\ell_1}|_x \ne T_{\ell_2}|_x$, and
    
  \item \ilabel{l4} if $\ell_1$ and $\ell_2$ are any two contact lines
    through $x$, then $\ell_1 \cap \ell_2 = \{x\}$.
  \end{enumerate}
\end{thm}

The smoothness of the base of the cone guarantees that much of the
theory developed by Hwang and Mok for uniruled varieties can be
applied to the contact setup. We refer to \cite{Hwa00} for an overview
and mention two examples.

\subsubsection{Stability of the tangent bundle}

It has been conjectured for a long time that complex contact manifolds
$X$ with $b_2(X)=1$ always carry a Kähler-Einstein metric. In
particular, it is conjectured that the tangent bundle of these
manifolds is stable. Using methods introduced by Hwang, stability
follows as an immediate consequence of Theorem~\ref{thm:main}.

\begin{cor}\label{cor:stability}
  Let $X$ be a complex-projective contact manifold with $b_2(X)=1$.
  Then the tangent bundle $T_X$ is stable.
\end{cor}

\subsubsection{Continuation of analytic morphisms}

The following corollary asserts that a contact manifold is determined
in a strong sense by the tangent directions to contact lines. The
analogous result for homogeneous manifolds appears in the work of
Yamaguchi.

\begin{cor}\label{cor:extension}
  Let $X$ be a complex-projective contact manifold and $X'$ be an
  arbitrary Fano manifold. Assume that $b_2(X) = b_2(X')=1$ and choose
  a dominating family of rational curves of minimal degree on
  $\mathcal H \subset \RatCurves(X')$. Assume further that there exist
  analytic open subsets $U \subset X$, $U' \subset X'$ and a
  biholomorphic morphism $\phi: U \to U'$ such that the tangent map
  $T\phi$ maps tangents of contact lines to tangents of curves coming
  from $\mathcal H$, and vice versa. Then $\phi$ extends to a
  biholomorphic map $\phi: X \to X'$.
\end{cor}

\begin{question}
  What would be the analogous statement in Riemannian geometry?
\end{question}

\subsection{Outline of this paper}

Property \iref{l1} of Theorem~\ref{thm:main} is known from previous
works ---see Fact~\ref{fact:perp} below. After a review of known facts
in chapter~\ref{chapt:facts}, properties \iref{l2}--\iref{l4} are
shown one by one in chapters~\ref{chapt:irred}--\ref{chapt:banana},
respectively. With these results at hand, the proofs of
Theorem~\ref{thm:main} and Corollaries~\ref{cor:stability} and
\ref{cor:extension}, which we give in
chapter~\ref{chapt:mainthmproof}, are very short.

The main difficulty in this paper is the proof of property \iref{l3},
which is done by a detailed analysis of the restriction of the tangent
bundle $T_X$ to pairs of contact lines that intersect tangentially.
The proof relies on a number of facts on jet bundles and on
deformation spaces of morphisms between polarized varieties for which
the author could not find any reference. These more general results
are gathered in the two appendices.

\subsection{Acknowledgements}

The main ideas for this paper were perceived while the author visited
the Korea Institute for Advanced Study in 2002. Details were worked
out during a visit to the University of Washington, Seattle, and while
the author was Professeur Invite at the Université Louis Pasteur in
Strasbourg. The author wishes to thank his hosts, Olivier Debarre,
Jun-Muk Hwang and Sándor Kovács for the invitations, and for many
discussions.

\section{Known Facts}\label{chapt:facts}

The proof of Theorem~\ref{thm:main} relies on a number of known facts
scattered throughout the literature.  For the reader's convenience, we
have gathered these results here.  Full proofs were included where
appropriate.

\subsection{Jet bundles on contact manifolds}\label{sec:jetoncontact}

The O'Neill tensor yields an identification $F \cong F^\vee\otimes L$.
If we dualize the contact sequence~\eqref{eq:contact-sequence} and
twist by $L$, we obtain a sequence,
\begin{equation}\label{eq:A}
  \begin{CD}
    0 @>>> \O_X @>>> \Omega^1_X \otimes L @>>> \underbrace{F}_{\cong
      F^\vee \otimes L} @>>> 0,
  \end{CD}
\end{equation}
which we would now like to compare to the dual of the first
jet sequence of $L$ ---see Appendix~\ref{app:jetsequence} for more
information on jets and the first jet-sequence.

By \cite[Thm.~2.1]{Leb95}, there exists a canonical symplectic form on
the $\mathbb C^*$-principal bundle associated with $L$ which gives
rise to an identification $\Jet^1(L) \cong \Jet^1(L)^\vee\otimes L$.
Thus, if we dualize the jet sequence and twist by $L$, we obtain a
sequence
\begin{equation}\label{eq:B}
  \begin{CD}
    0 @>>> \O_X @>>> \underbrace{\Jet^1(L)}_{\cong \Jet^1(L)^\vee
      \otimes L} @>>> T_X @>>> 0
  \end{CD}
\end{equation}
It is known that sequence~\eqref{eq:A} is a sub-sequence of
\eqref{eq:B}.

\begin{fact}[{\cite[p.~426]{Leb95}}]\label{fact:jetsoncontact}
  There exists a commutative diagram with exact rows and columns
  $$
  \xymatrix{ 
    & & & 0 \ar[d] & 0 \ar[d] \\
    {\eqref{eq:A}} & 0 \ar[r] & {\O_X} \ar@{=}[d] \ar[r] & 
    {\Omega^1_X \otimes L} \ar[r] \ar[d] & {F} \ar[r] \ar[d] & 0 \\
    {\eqref{eq:B}} & 0 \ar[r] & {\O_X} \ar[r] & {\Jet^1(L)} \ar[r]
    \ar[d] & {T_X} \ar[r] \ar[d] & 0 \\
    & & & L \ar@{=}[r] \ar[d] & L \ar[d] \\
    & & & 0 & 0 \\
  }  
  $$
  where the middle column is the first jet sequence for $L$ and the
  right column is the sequence~\eqref{eq:contact-sequence} of
  page~\pageref{eq:contact-sequence} that defines the contact
  structure.
\end{fact}

\subsection{Contact Lines}\label{sect:contLines}

It is conjectured that a projective contact manifold $X$ with
$b_2(X)=1$ is homogeneous. This is known to be equivalent to
conjecture that there exists an embedding $X \to \P^N$ that maps
contact lines to lines in $\P^N$. While we cannot presently prove
these conjectures, it has already been shown in the first part
\cite{Kebekus01} of this work that a contact lines through a general
point share many features with lines in $\P^N$. Some of the following
results will be strengthened in Chapter~\ref{sect:31}.

\begin{fact}[{\cite[Rem.~3.3]{Kebekus01}}]\label{fact:Fintl}
  Let $\ell$ be a contact line. Then $\ell$ is $F$-integral. In other
  words, if $x \in \ell$ is a smooth point, then $T_\ell|_x \subset
  F|_y$.
\end{fact}

\begin{fact}\label{fact:perp}
  Let $x \in X$ be a general point and $\ell \subset X$ a contact line
  that contains $x$. Then $\ell$ is smooth. The splitting types of
  the restricted vector bundles $F|_\ell$ and $T_X|_\ell$ are:
  \begin{align*}
    T_X|_\ell & \cong \O_\ell(2) \oplus \O_\ell(1)^{\oplus n-1} \oplus \O_\ell^{\oplus n+1} \\
    F|_\ell   & \cong \underbrace{\O_\ell(2) \oplus \O_\ell(1)^{\oplus n-1} \oplus
      \O_\ell^{\oplus n-1}}_{=: F|_\ell^{\geq 0}}\oplus \, \O_\ell(-1)
  \end{align*}
  For all points $y \in \ell$, the vector space $F|_\ell^{\geq 0} |_y$
  and the tangent space $T_\ell|_y$ are perpendicular with respect to
  the O'Neill tensor $N$: $F|_\ell^{\geq 0} |_y = T_\ell|_y^\perp$.
\end{fact}
\begin{proof}
  The fact that $\ell$ is smooth was shown in
  \cite[Prop.~3.3]{Kebekus01}.  The splitting type of $T_X|_\ell$ is
  given by \cite[Lem.~3.5]{Kebekus01}. To find the splitting type of
  $F|_\ell$, recall that the contact structure yields an
  identification $F \cong F^\vee \otimes L$. Since $L|_\ell \cong
  \O_\ell(1)$, we can therefore find positive numbers $a_i$ and write
  $$
  F|_\ell \cong \bigoplus^n_{i=1} (\O_\ell(a_i) \oplus \O_\ell(1-a_i)) 
  $$
  The precise splitting type then follows from the splitting type
  of $T_X|_\ell$ and from Fact~\ref{fact:Fintl} above.
  
  The simple observation that every map $\O_\ell(2) \cong T_\ell \to
  L|_\ell \cong \O_\ell(1)$ is necessarily zero yields the fact that
  $F|_\ell^{\geq 0} |_y$ and $T_\ell|_y$ are perpendicular with
  respect to the O'Neill tensor $N$.
\end{proof}

\begin{fact}\label{fact:3}
  Let $x \in X$ be a general point, $\ell \subset X$ a contact line
  that contains $x$ and $y \in \ell$ any point. If $s \in H^0(\ell,
  T_X|_\ell)$ is a section such that $s(y) \in F|_y$, then $s$ is
  contained in $H^0(\ell, F|_\ell)$ if and only if $T_\ell|_y$ and
  $s(y)$ are orthogonal with respect to the O'Neill-tensor $N$.
  
  In particular, we have that if $s(y) \in F|_\ell^{\geq 0}$, then $s
  \in H^0(\ell, F|_\ell^{\geq 0})$.
\end{fact}
\begin{proof}
  Let $f : \P^1 \to X$ be a parameterization of $\ell$. We know from
  \cite[Thms.~II.3.11.5 and II.2.8]{K96} that the space $\Hom(\P^1,X)$
  is smooth at $f$. Consequence: we can find an embedded unit disc
  $\Delta \subset \Hom(\P^1,X)$, centered about $f$ such that $s \in
  T_\Delta|_f$ holds ---see Fact~\ref{fact:tgt_vect_to_Hom} on page
  \pageref{fact:tgt_vect_to_Hom} for a brief explanation of the
  tangent space to $\Hom(\P^1,X)$. In this situation we can apply
  \cite[Prop.~3.1]{Kebekus01} to the family $\Delta$, and the claim is
  shown.
\end{proof}

\subsection{Dubbies}

In Section~\ref{sec:proofOfThm1} we will show that no two contact
lines through a general point share a common tangent direction at $x$.
For this, we will argue by contradiction and assume that $X$ is
covered by pairs of contact lines which intersect tangentially in at
least one point. Such a pair is always dominated by a pair of smooth
rational curves that intersect in one point with multiplicity exactly
2. These particularly simple pairs were called ``dubbies'' and
extensively studied in \cite[Sect.~3]{KK02}.

\begin{defn}\label{def:dubby}
  A \emph{dubby} is a reduced, reducible curve, isomorphic to the
  union of a line and a smooth conic in $\P^2$ intersecting
  tangentially in a single point.
\end{defn}

\begin{rem}\label{rem:dubbyconic}
  Because a dubby $\ell$ is a plane cubic, we have $h^1(\ell,
  \O_\ell)=1$.
\end{rem}

\subsubsection{Low degree line bundles on dubbies}

It is the key observation in \cite{KK02} that an ample line bundle on
a dubby always gives a canonical identification of its two irreducible
components. In the setup of section~\ref{sec:proofOfThm1}, where
dubbies are composed of contact lines, the identification is quite
apparent so that we do not need to refer to the complicated general
construction of \cite[Sect.~3.2]{KK02}.

\begin{prop}\label{prop:mirrors}
  Let $\ell = \ell_1 \cup \ell_2$ be a dubby and $H \in \Pic(\ell)$ be
  a line bundle whose restriction to both $\ell_1$ and $\ell_2$ is of
  degree one. Then there exists a unique isomorphism $\gamma: \ell \to
  \P^1$ such that
  \begin{enumerate}
  \item the restriction $\gamma|_{\ell_i} : \ell_i \to \P^1$ to any
    component is isomorphic and
  \item a pair of smooth points $y_1 \in \ell_1$ and $y_2 \in \ell_2$
    forms a divisor for $H$ if and only if $\gamma(y_1)= \gamma(y_2)$.
  \end{enumerate}
  In particular, we have that $h^0(\ell, H)=h^0(\P^1,
  \O_{\P^1}(1))=2$.
\end{prop}
\begin{proof}
  Consider the restriction morphisms
  $$
  r_i : H^0(\ell, H) \to H^0(\ell_i, H|_{\ell_i}) \simeq
  H^0(\P^1,\O_{\P^1}(1)).
  $$
  We claim that the morphism $r_i$ is an isomorphism for all
  $i\in\{1,2\}$.  The rôles of $r_1$ and $r_2$ are symmetric, so it is
  enough to prove the claim for $r_1$.  First note that $h^0(\ell, H)
  \geq 2$ by \cite[Lem.~3.2]{KK02}. It is then sufficient to prove
  that $r_1$ is injective.  Let $s \in \ker(r_1) \subset H^0(\ell,
  H)$. In order to show that $s=0$ it is enough to show that $r_2(s) =
  0$. Notice that $r_2(s)$ is a section in $H^0(\ell_2, H|_{\ell_2})$
  that vanishes on the scheme-theoretic intersection $\ell_1 \cap
  \ell_2$.  The length of this intersection is two and any non-zero
  section in $H^0(\ell_2, H|_{\ell_2}) \simeq H^0(\P^1,\O_{\P^1}(1))$
  has a unique zero of order one, hence $r_2(s)$ must be zero, and so
  $r_i$ is indeed an isomorphism for all $i\in\{1,2\}$.
  
  This implies that $H$ is generated by global sections and gives a
  morphism $\gamma: \ell \to \P^1$, whose restriction
  $\gamma|_{\ell_i}$ to any of the two components is an isomorphism.
  Property \iref{l2} follows by construction.
\end{proof}
\begin{notation}
  We call a pair of points $(y_1, y_2)$ as in
  Proposition~\ref{prop:mirrors} ``mirror points with respect to
  $H$''.
\end{notation}

\begin{cor}\label{cor:actiononpic}
  Let $\ell = \ell_1 \cup \ell_2$ be a dubby and $\Pic^{(1,1)}(\ell)$
  be the component of the Picard-group that represents line bundles
  whose restriction to both $\ell_1$ and $\ell_2$ is of degree one.
  Then the natural action of the automorphism group $\Aut(\ell)$ on
  $\Pic^{(1,1)}(\ell)$ is transitive.
\end{cor}
\begin{proof}
  Consider the open set $\Omega = \ell_2 \setminus (\ell_1 \cap
  \ell_2)$. By Proposition~\ref{prop:mirrors} it suffices to show that
  there exists a group $G \subset \Aut(\ell)$ that fixes $\ell_1$
  pointwise and acts transitively on $\Omega$.
  
  For this, define a group action on the disjoint union $\ell_1
  \coprod \ell_2$ as follows. Let $G \subset \Aut(\ell_2)$, $G \cong
  \mathbb C$ be the isotropy group of the scheme-theoretic
  intersection $\ell_1 \cap \ell_2 \subset \ell_2$. Let $G$ act
  trivially on $\ell_1$. It is clear that $G$ acts freely on $\Omega$.
  By construction, $G$ acts trivially on the scheme-theoretic
  intersection $\ell_1 \cap \ell_2$ so that the actions on $\ell_1$
  and $\ell_2$ glue to give a global action on $\ell$.
\end{proof}

\begin{cor}\label{cor:autlh}
  Let $\ell$ and $H$ be as in Proposition~\ref{prop:mirrors} above and
  let
  $$
  \Aut(\ell, H) := \{g \in \Aut(\ell) \,|\, g^*(H)\cong H\}
  $$
  be the subgroup of automorphisms that respect the line bundle
  $H$. If $y \in \ell$ is any smooth point, then there exists a vector
  field, i.e., a section of the tangent sheaf
  $$
  s \in T_{\Aut(\ell, H)}|_{\rm Id} \subset H^0(\ell, T_\ell)
  $$
  that does not vanish at $y$.
\end{cor}
\begin{proof}
  Let $\sigma = \ell_1 \cap \ell_2$ be the (reduced) singular point,
  let $\eta: \ell_1 \coprod \ell_2 \to \ell$ be the normalization and
  consider the natural action of $\C$ on $\P^1$ that fixes the image
  point $\gamma(\sigma) \in \P^1$. Use the isomorphisms
  $\gamma|_{\ell_1}$ and $\gamma|_{\ell_2}$ to define a $\C$-action on
  $\ell_1 \coprod \ell_2$. As before, observe that this action acts
  trivially on the scheme-theoretic preimage
  $$
  \eta^{-1}(\ell_1 \cap \ell_2).
  $$
  The $\C$-action on $\ell_1 \coprod \ell_2$ therefore descends to a
  $\C$-action on $\ell$. To see that the associated vector field does
  not vanish on $y$, it suffices to note that the singular point
  $\sigma$ is the only $\C$-fixed point on $\ell$. 

  Because the action preserves $\gamma$-fibers, it follows from
  Proposition~\ref{prop:mirrors} that $\C$ acts via a morphism
  $$
  \C \to \Aut(\ell, H).
  $$
\end{proof}

In section~\ref{sec:proofOfThm1} we will need to consider line bundles
of degree $(2,2)$. The following remark will come handy.

\begin{lem}\label{lem:bideg22}
  Let $\ell = \ell_1 \cup \ell_2$ be a dubby and let $H \in
  \Pic(\ell)$ be a line bundle whose restriction to both $\ell_1$ and
  $\ell_2$ is of degree 2. For $i \in \{1,2\}$ there exist sections
  $s_i \in H^0(\ell, H)$ which vanish identically on $\ell_i$ but not
  on the other component.
\end{lem}
\begin{proof}
  By \cite[Lem.~3.2]{KK02}, we have $h^0(\ell, H) \geq 4$. Thus, the
  restriction map $H^0(\ell, H) \to H^0(\ell_i, H|_{\ell_i}) \cong
  \C^3$ has a non-trivial kernel.
\end{proof}

\subsubsection{Vector bundles on dubbies}

Dubbies are in many ways similar to elliptic curves. While $H^1(\ell,
\O_\ell)$ does not vanish, the higher cohomology groups of ample
vector bundles are trivial.

\begin{lem}\label{lem:vanishingForDubbies}
  Let $\mathcal E$ be a vector bundle on $\ell$ whose restrictions to
  both $\ell_1$ and $\ell_2$ is ample. Then $H^1(\ell, \mathcal E)=0$.
\end{lem}
\begin{proof}
  Let $\overline \sigma := \ell_1 \cap \ell_2 \subset \ell$ be the
  scheme-theoretic intersection, which is a zero-dimensional subscheme
  of length two. Now consider the normalization $\eta: \ell_1 \coprod
  \ell_2 \to \ell$ and the associated natural sequence
  \begin{equation}\label{eq:conductorseq}
    \begin{CD}
      0 @>>> \mathcal E @>>> \eta_*(\eta^* \mathcal E) @>{\alpha}>>
      \mathcal E|_{\overline \sigma} @>>> 0
    \end{CD}
  \end{equation}
  where $\alpha$ is defined on the level of pre-sheaves as follows.
  Assume we are given an open neighborhood $U$ of the singular point
  $\sigma \in \ell$. By definition of $\eta_*(\eta^* \mathcal E)$, to
  give a section $s \in \eta_*(\eta^* \mathcal E)(U)$ it is equivalent
  to give two sections $s_1 \in (\mathcal E|_{\ell_1})(U \cap \ell_1)$
  and $s_2 \in (\mathcal E|_{\ell_2})(U \cap \ell_2)$. If
  $$
  r_i : (\mathcal E|_{\ell_i})(U \cap \ell_i) \to \mathcal E|_{\overline
    \sigma}
  $$
  are the natural restriction morphisms, then we write $\alpha$ as
  $$
  \alpha(s) = r_1(s_1)-r_2(s_2).
  $$
  A section of the long homology sequence associated
  with~\eqref{eq:conductorseq} reads
  $$
  \begin{CD}
    H^0(\ell, \eta_* \eta^* \mathcal E) @>{\beta}>> H^0(\overline
    \sigma, \mathcal E|_{\overline \sigma}) @>>> H^1(\ell, \mathcal E) @>>>
    H^1(\ell, \eta_* \eta^* \mathcal E),
  \end{CD}
  $$
  where $\beta$ is again the difference of the restriction
  morphisms. We have that
  \begin{align*}
    H^0(\ell, \eta_* \eta^* \mathcal E) &= H^0(\ell_1, \mathcal
    E|_{\ell_1}) \oplus H^0(\ell_2, \mathcal E|_{\ell_2}) \\
    H^1(\ell, \eta_* \eta^* \mathcal E) &= H^1(\ell_1, \mathcal
    E|_{\ell_1}) \oplus H^1(\ell_2, \mathcal E|_{\ell_2}) = \{0\}
  \end{align*}
  and it remains to show that $\beta$ is surjective. That, however,
  follows from the fact that $\mathcal E|_{\ell_i}$ is an ample bundle
  on $\P^1$ that generates 1-jets so that even the single restriction
  $$
  r_1 : H^0(\ell_1, \mathcal E|_{\ell_1}) \to H^0(\overline \sigma,
  \mathcal E|_{\overline \sigma})
  $$
  alone is surjective.
\end{proof}

\section{Irreducibility}\label{chapt:irred}

As a first step towards the proof of Theorem~\ref{thm:main}, we
show the irreducibility of the space of contact lines through a
general point.

\begin{thm}\label{thm:irreducibility}
  If $x \in X$ is a general point, then the subset $H_x \subset H$ of
  contact lines through $x$ is irreducible. In particular,
  $\locus(H_x)$ is irreducible.
\end{thm}

The proof of Theorem~\ref{thm:irreducibility}, which is given in
Section~\ref{sect:32} below, requires a strengthening of
Fact~\ref{fact:perp}, which we give in the following section.

\subsection{Contact lines with special splitting type}\label{sect:31}

We adopt the notation of \cite[Chapt.~1.2]{Hwa00} and call a contact
line $\ell \subset X$ ``standard'' if
$$
\eta^*(T_X) \cong \O_{\P^1}(2) \oplus \O_{\P^1}(1)^{\oplus n-1}
\oplus \O_{\P^1}^{\oplus n+1},
$$
where $\eta: \P^1 \to \ell$ is the normalization. It is known that
the set of standard curves is Zariski-open in $H$, see again
\cite[Chapt.~1.2]{Hwa00}. We can therefore consider the subvariety
$$
H' := \{ \ell \in H \,|\, \text{ $\ell$ not standard} \}
$$
The proof of Theorem~\ref{thm:irreducibility} is based on the
observation that there is only a small set in $X$ whose points are not
contained in a standard contact line. For a proper formulation, set 
$$
D := \locus(H') = \bigcup_{\ell \in H'} \ell.
$$
If follows immediately from Fact~\ref{fact:perp} that $D$ is a proper
subset of $X$.

\begin{prop}\label{prop:D0}
  If $D^0 \subset D$ is any irreducible component with $\codim_X D^0
  =1$, $x \in D^0$ a general point, and $H^0_x \subset H_x$ any
  irreducible component, then there exists a curve $\ell \in H^0_x$
  which is not contained in $D$ and therefore free.
\end{prop}

The proof of Proposition~\ref{prop:D0} is a variation of the
argumentation in \cite[Chapt.~4]{Kebekus01}. While we work here in a
more delicate setup, moving only components of $\locus(H_x)$ along the
divisor $D^0$, parts of the proof were taken almost verbatim from
\cite{Kebekus01}.

\subsubsection*{Proof of Proposition~\ref{prop:D0}, Step 1: Setup} 
Assume to the contrary, i.e., assume that there exists a divisor $D^0
\subset D$ such that for a general point $x \in D^0$ there exists a
component of $H_x$ whose associated curves are all contained in $D^0$.
Since by \cite[Prop.~4.1]{Kebekus01} for all $y \in X$, the space
$H_y$ is of pure dimension $n-1$, we can find a closed, proper
subvariety $H^0 \subset H$ with $\locus(H^0) = D^0$ such that for all
points $y \in D^0$,
$$
H^0_y = \{ \ell \in H^0 \,|\, y \in \ell \}
$$
is the union of irreducible components of $H_y$. In particular, we
have that for all $y \in D^0$, $\dim \locus (H^0_y) = n$.

\subsubsection*{Proof of Proposition~\ref{prop:D0}, Step 2:
  Incidence variety}
In analogy to \cite[Notation~4.2]{Kebekus01}, define the incidence
variety
$$
V := \{ (x',x'') \in D^0\times D^0 \,|\, x'' \in \locus(H^0_{x'}) \}
\subset D^0 \times X.
$$
Let $\pi_1, \pi_2 : V \to D^0$ be the natural projections. We have
seen in Step~1 above that for every point $y \in D^0$, $\pi^{-1}_1(y)$
is a subscheme of $X$ of dimension $\dim \pi^{-1}_1(y) = n$. In
particular, $V$ is a well-defined family of cycles in $X$ in the sense
of \cite[Chapt.~I.3.10]{K96}. The universal property of the
Chow-variety therefore yields a map
$$
\Phi : D^0 \to \Chow(X).
$$
Since $\dim \locus(H^0_y) = n < \dim D^0$, this map is not
constant. On the other hand, since $D^0 \subset X$ is ample, the
Lefschetz hyperplane section theorem \cite[Thm.~2.3.1]{BS95} asserts
that $b_2(D^0)=1$. As a consequence, we obtain that the map $\Phi$ is
finite: for any given point $y \in D^0$ there are at most finitely
many points $y_i$ such that $\locus(H^0_y) = \locus(H^0_{y_i})$. In
analogy to \cite[Lemma~4.3]{Kebekus01} we conclude the following.

\begin{lem}\label{lem:gamma}
  Let $\Delta$ be a unit disc and $\gamma: \Delta \to D^0$ an
  embedding. Then there exists a Euclidean open set $V^0 \subset
  \pi^{-1}_1 (\gamma(\Delta))$ such that $\pi_2(V^0) \subset X$ is a
  submanifold of dimension
  $$
  \dim \pi_2(V^0) = n+1.
  $$
  In particular, $\pi_2(V^0)$ is not $F$-integral. \qed
\end{lem}

\subsubsection*{Proof of Proposition~\ref{prop:D0}, Step 3:
  conclusion}

We shall now produce a map $\gamma:\Delta \to D^0$ to which
Lemma~\ref{lem:gamma} can be applied. For that, recall that $D^0$
cannot be $F$-integral. Thus, if $y \in D^0$ is a general smooth point
of $D^0$, then
$$
F_{D^0,y} := F|_y \cap T_{D^0}|_y
$$
is a proper hyperplane in $F|_y$, and the set $F_{D^0,y}^\perp$ of
tangent vector that are orthogonal to $F_{D^0,y}$ with respect to the
O'Neill-tensor is a line that is contained in $F_{D^0,y}$. The
$F_{D^0,y}$ give a (singular) 1-dimensional foliation on $D^0$ which
is regular in a neighborhood of the general point $y$. Let $\gamma:
\Delta \to D^0$ be an embedding of the unit disk that is an integral
curve of this foliation, i.e., a curve such that for all points $y'
\in \gamma(\Delta)$ we have that
\begin{equation}\label{form:perp_assert}
  T_{\gamma(\Delta)}|_{y'} = F_{D^0,y'}^\perp
\end{equation}
Now let $\mathcal H \subset (\Hom_{bir}(\P^1,X))_{\red}$ be the family
of generically injective morphisms parameterizing the curves
associated with $H^0$. Fix a point $0 \in \P^1$ and set
$$
\mathcal H_{\Delta} := \{f\in \mathcal H\ |\ f(0)\in
\gamma(\Delta)\}.
$$
If $\mu : \mathcal H_\Delta \times \P^1 \to X$ is the universal
morphism, then it follows by construction that
$$
\mu (\mathcal H_\Delta\times \P^1) = 
\pi_2(\pi_1^{-1}(\gamma(\Delta)))\supset\pi_2(V^0),
$$
where $V^0$ comes from Lemma~\ref{lem:gamma}. In particular, since
$\pi_2(V^0)$ is not $F$-integral, there exists a smooth point
$(f,p)\in \mathcal H_\Delta\times \P^1$ with $f(p) \in \pi_2(V^0)$ and
there exists a tangent vector $\vec w\in T_{\mathcal H_\Delta\times
  \P^1}|_{(f,p)}$ such that the image of the tangent map is not in
$F$:
$$
T{\mu_\Delta}(\vec w) \not \in  \mu^*(F).
$$
Decompose $\vec w = \vec w'+\vec w''$, where $\vec w\in
T_{\P^1}|_p$ and $\vec w''\in T_{\mathcal H_\Delta}|_f$. Then, since
$f(\P^1)$ is $F$-integral, it follows that $T\mu(\vec w')\in \mu^*(F)$
and therefore
\begin{equation}\label{eq:wprimeprime}
  T\mu(\vec w'')\not \in \mu^*(F).    
\end{equation}
As a next step, since $\mathcal H_\Delta$ is smooth at $f$, we can
choose an immersion
$$
\begin{array}{rrcl}
  \beta : & \Delta & \to     & \mathcal H_\Delta \\
  & t      & \mapsto & \beta_t
\end{array}
$$
such that $\beta_0=f$ and such that
$$
T\beta\left(\left. \frac{\partial}{\partial t}\right|_{t=0}\right) 
= \vec w''.
$$
In particular, if $s \in H^0(\P^1, f^*(T_X))$ is the section
associated with $\vec w'' = T\beta(\frac{\partial}{\partial
  t}|_{t=0})$, and $s' := f^*(\theta)(s)\in H^0(\P^1,f^*(L))$, then
the following holds:
\begin{enumerate}
\item it follows from (\ref{eq:wprimeprime}) and from
  \cite[Prop.~II.3.4]{K96} that $s'$ is not identically zero.
  
\item at $0 \in \P^1$, the section $s$ satisfies $s(0) \in
  f^*(T_{\gamma(\Delta)})\subset f^*(F)$. In particular, $s'(0) = 0$.
    
\item If $z$ is a local coordinate on $\P^1$ about $0$, then it
  follows from \eqref{form:perp_assert} that $\frac{\partial}{\partial
    z}|_0 \in f^*(F)$ and $s(0) \in f^*(F')$ are perpendicular with
  respect to the non-degenerate form $N$.
\end{enumerate}
  
Items \iref{l2} and \iref{l3} ensure that we can apply
\cite[Prop.~3.1]{Kebekus01} to the family $\beta_t$.  Since the
section $s'$ does not vanish completely, the proposition states that
$s'$ has a zero of order at least two at $0$. But $s'$ is an element
of $H^0(\P^1,f^*(L))$, and $f^*(L)$ is a line bundle of degree one. We
have thus reached a contradiction, and the proof of
Proposition~\ref{prop:D0} is finished. \qed

\subsection{Proof of Theorem~\ref{thm:irreducibility}}\label{sect:32}

Let $\pi :U \to H$ be the restriction of the universal $\P^1$-bundle
$\Univrc(X)$ to $H$ and let $\iota: U \to X$ be the universal
morphism. Consider the Stein-factorization of $\iota$.
$$
\xymatrix{
  U \ar@/^0.6cm/[rrrr]^{\iota = \alpha \circ \beta} 
  \ar[rr]^{\alpha}_{\text{connected fibers}} 
  \ar[d]_{\text{$\P^1$-bundle $\pi$}} & & 
  {X'} \ar[rr]^{\beta}_{\text{finite}} && X \\
  H }
$$
Let $T \subset X$ be the union of the following subvarieties of $X$:
\begin{itemize}
\item the components $D_i \subset D$ which have $\codim_X D_i \geq 2$,
  where $D \subset X$ is the subvariety defined in
  section~\ref{sect:31} above.
  
\item for every divisorial component $D_i \subset D$, the
  Zariski-closed set of points $y \in D_i$ for which there exists an
  irreducible component $H_y^0 \subset H_y$ such that none of the
  associated curves in $X$ are free
  
\item the image $\beta(X'_{\Sing})$ of the singular set of $X'$
\end{itemize}
It follows immediately from Proposition~\ref{prop:D0} that $\codim_X T
\geq 2$.

\begin{claim}\label{claim:smoothnessofbeta}
  The morphism $\beta$ is unbranched away from $T$, i.e., the
  restricted morphism 
  $$
  \beta|_{X\setminus T} : \beta^{-1}(X \setminus T) \to X \setminus T
  $$
  is smooth.
\end{claim}

\subsubsection*{Proof of Claim~\ref{claim:smoothnessofbeta}}
Let $y \in \beta^{-1}(X \setminus T)$ be any point. To show that
$\beta$ has maximal rang at $y$, it suffices to find a point $z \in
\alpha^{-1}(y)$ such that 
\begin{itemize}
\item $z$ is a smooth point of $U$ and such that
\item $\iota$ is smooth at $z$.
\end{itemize}
By \cite[Chapt. II, Thms. 1.7, 2.15 and Cor. 3.5.4]{K96}, both
requirements are satisfied if $\pi(z) \in H$ is a point that
corresponds to a free curve. The existence of a free curve in the
component $\pi(\alpha^{-1}(y))$, however, is guaranteed by choice of
$T$. \qed

\subsubsection*{Application of Claim~\ref{claim:smoothnessofbeta}}
Since $X$ is Fano, it is simply connected. Because $T \subset X$ is
not a divisor, its complement $X\setminus T$ is also simply connected.
Claim~\ref{claim:smoothnessofbeta} therefore implies that $X'$ is
either reducible, or that the general $\beta$-fiber is a single point.
But $X'$ is irreducible by construction, and it follows that the
general fiber of $\iota$ must be connected. By Seidenberg's classical
theorem \cite[Thm.~1.7.1]{BS95}, the general fiber $\iota^{-1}(x)$ is
then irreducible, and so is its image
$$
H_x = \pi(\iota^{-1}(x)).
$$
This ends the proof of Theorem~\ref{thm:irreducibility}. \qed

\section{Contact lines sharing a common tangent direction}
\label{sec:proofOfThm1}

The aim of the present section~\ref{sec:proofOfThm1} is to give a
proof of part \iref{l3} of Theorem~\ref{thm:main}. More precisely, we
show the following.

\begin{thm}\label{thm:main2and3}
  If $x \in X$ is a general point, then all contact lines through $x$
  are smooth, and no two of them share a common tangent at $x$.
\end{thm}

The proof is at its core a repeat performance of
\cite[Sect.~4.1]{KK02} where the global assumptions of
\cite[Thm.~1.3]{KK02} are replaced by a careful study of the
restriction of the tangent bundle $T_X$ to a pair of rational curves
with non-transversal intersection.

\subsection{Setup}\label{sect:4setup}

We will argue by contradiction and assume throughout the rest of this
section to the contrary. More precisely, we stick to the following.
\begin{assumption}\label{ass:thm11}
  Assume that at for a general point $x \in X$, we can find a pair
  $\ell' = \ell'_1 \cup \ell'_2 \subset X$ of distinct contact lines
  $\ell'_i \in H$ that intersect tangentially at $x$.
\end{assumption}
The pair $\ell'$ is then dominated by a dubby $\ell = \ell_1 \cup
\ell_2$ whose singular point $\sigma = \ell_1 \cap \ell_2$ maps to
$x$. For the remainder of this section we fix a generically injective
morphism $f : \ell \to \ell'$ such that $f(\sigma)=x$. We also fix the
line bundle $H := f^*(L) \in \Pic^{(1,1)}(\ell)$.

\subsection{Proof of Theorem~\ref{thm:main2and3}}
The Assumption~\ref{ass:thm11} implies that for a fixed point $x$,
there is a positive dimensional family of pairs curves which contain
$x$ and have a point of non-transversal intersection. Loosely
speaking, we will move the point of intersection to obtain a
positive-dimensional family of dubbies that all contain the point $x$.

To formulate more precisely, consider the quasi-projective reduced
subvariety
$$
\mathcal H \subset (\Hom(\ell,X))_{\red}
$$
of morphisms $g \in \Hom(\ell,X)$ such that $g^*(L) \cong H$. Note
that such a morphism will always be generically injective on each
irreducible component of $\ell$. Consider the diagram
$$
\xymatrix{
  {\mathcal H \times \ell} \ar[d]_{\text{projection $\pi$}}
  \ar[rrr]^{\mu}_{\text{universal morphism}} & & &  {X} \\
  {\mathcal H} }
$$
and conclude from Corollary~\ref{cor:actiononpic} that the
restricted universal morphism $\mu|_{\mathcal H \times \{\sigma\}}$ is
dominant.  By general choice of $f$, there exists a unique
positive-dimensional irreducible component
$$
\mathcal H_x \subset \pi ( \mu^{-1}(x) )
$$
which contains $f$ and which is smooth at $f$. It is clear that for
a general point $g \in \mathcal H_x$, the point $x$ is a smooth point
of the pair of curves $g(\ell)$. This implies the following
decomposition lemma.

\begin{lem}\label{lem:decomp}
  The preimage of $x$ decomposes as
  $$
  \mu^{-1}(x) \cap \pi^{-1}(\mathcal H_x) = \tau_1 \cup \bigcup_{i=1}^N \eta_i,  
  $$
  where $\tau_1 \subset \mathcal H_x \times \ell$ is a section that
  intersects $\mathcal H_x \times \{\sigma \}$ over $f$, but is not
  contained in $\mathcal H_x \times \{\sigma \}$, and where the
  $\eta_i$ are finitely many lower-dimensional components, $\dim
  \eta_i < \dim H_x$.
\end{lem}

\begin{proof}[Proof of Lemma~\ref{lem:decomp}]
  Since all curves in $X$ that are associated with points of $\mathcal
  H_x$ contain $x$, it is clear that there exists a component $\tau_1
  \subset \mu^{-1}(x) \cap \pi^{-1}(\mathcal H_x)$ that surjects onto
  $\mathcal H_x$.
  
  We have seen above, that for $g \in \mathcal H_x$ general, $x$ is a
  smooth point of the pair of curves $g(\ell)$, i.e.~that the
  scheme-theoretic intersection $\mu^{-1}(x) \cap \pi^{-1}(g)$ is a
  single closed point that is not equal to $\sigma$. Since
  $\mu^{-1}(x) \cap \pi^{-1}(g)$ is necessarily discrete for all $g
  \in \mathcal H_x$, it follows that $\tau_1$ is a section that is not
  contained in $\mathcal H_x \times \{\sigma \}$. It follows further
  that no other component $\eta_i$ of $\mu^{-1}(x) \cap
  \pi^{-1}(\mathcal H_x)$ dominates $\mathcal H_x$. In particular,
  $\dim \eta_i < \dim \tau_1$.
  
  To see that $(f,\sigma) \in \tau_1$, we first note that
  $f(\sigma)=0$, so that $(f,\sigma)$ is contained in the preimage,
  $(f,\sigma) \in \mu^{-1}(x) \cap \pi^{-1}(\mathcal H_x)$. On the
  other hand, Fact~\ref{fact:perp} of page~\pageref{fact:perp} asserts
  that both $f(\ell_1)$ and $f(\ell_2)$ are smooth so that $\sigma =
  f^{-1}(x)$ and $(f,\sigma) = \mu^{-1}(x) \cap \pi^{-1}(f)$. This
  ends the proof.
\end{proof}

After renaming $\ell_1$ and $\ell_2$, if necessary, we assume without
loss of generality that $\tau_1 \subset \mathcal H_x \times \ell_1$.
By Proposition~\ref{prop:mirrors}, the line bundle $H \in \Pic(\ell)$
yields an identification morphism $\gamma: \ell \to \P^1$. Let
$$
(Id \times \gamma) : \mathcal H_x \times \ell \to
\mathcal H_x \times \P^1
$$
be the associated morphism of bundles and consider the mirror
section
$$
\tau_2 := ((Id \times \gamma)|_{\mathcal H_x \times \ell_2})^{-1} (Id
\times \gamma)(\tau_1).
$$

\begin{claim}\label{claim:tau2}
  The universal morphism $\mu$ contracts $\tau_2$ to a point:
  $\mu(\tau_2)= (*)$.
\end{claim}

\begin{proof}[Proof of Claim~\ref{claim:tau2}]
  The proof of Claim~\ref{claim:tau2} makes use of
  Proposition~\ref{prop:section2} of page~\pageref{prop:section2}
  which is shown independently in
  sections~\ref{sec:subb}--\ref{sec:secinsubb} below.
  
  For the proof, we pick a general smooth point $z \in \tau_2$, and an
  arbitrary tangent vector $\vec v \in T_{\tau_2}|_z$. It suffices to
  show that $\vec v$ is mapped to zero,
  $$
  T\mu(\vec v) = 0 \in \mu^*(T_X)|_z.
  $$
  
  Since $\tau_2$ is a section over $\mathcal H_x$, and since $\mathcal
  H_x$ is smooth at $\pi(z)$, we can find a small embedded unit disc
  $\Delta \subset \mathcal H_x$ with coordinate $t$ such that
  $T\pi(\vec v) = \pi^*\left( \frac{\partial}{\partial t} \right)|_z$.
  For the remainder of the proof, it is convenient to introduce new
  bundle coordinates on the restricted bundle $\Delta \times \ell$. It
  follows from Corollary~\ref{cor:autlh} that, after perhaps shrinking
  $\Delta$, we can find a holomorphic map
  $$
  \alpha : \Delta \to \Aut^0(\ell, H)
  $$
  with associated coordinate change diagram
  $$
  \xymatrix{ {\Delta\times \ell} \ar[d]_{\text{projection $\pi_1$}}
    \ar[rr]^{\text{coord. change $\kappa$}}_{(t,y) \mapsto
      (t,\alpha(t)\cdot y)} \ar@/^0.7cm/[rrrrr]^{\tilde \mu := \mu
      \circ \kappa} && {\mathcal H \times \ell}
    \ar[rrr]^{\mu}_{\text{universal morphism}} & & &  {X} \\
    {\Delta} \ar[r]^(.4){\alpha} & {\Aut(\ell,H)}}
  $$
  such that $\kappa^{-1}(\tau_1) \cup \kappa^{-1}(\tau_2)$ is a
  fiber of the projection $\pi_2: \Delta \times \ell \to \ell$.
  
  Let $\tau'_i := \kappa^{-1}(\tau_i)$, $z' := \kappa^{-1}(z)$, and
  let $\vec v' \in T_{\tau'_2}|_{z'}$ be the preimage of $\vec v$,
  i.e.  the unique tangent vector that satisfies $T\kappa(\vec v') =
  \kappa^{-1}(\vec v)$. The new coordinates make it easy to write down
  an extension of the tangent vector $\vec v'$ to a global vector
  field, i.e. to a section $s \in H^0(\Delta \times \ell, T_{\Delta
    \times \ell})$ of the tangent sheaf. Indeed, if we use the product
  structure to decompose
  $$
  T_{\Delta \times \ell} \cong \pi_1^*(T_\Delta) \oplus \pi_2^*(T_\ell),
  $$
  then the ``horizontal vector field'' $s:=\pi_1^* \left(
    \frac{\partial}{\partial t} \right)$ will already satisfy $s(z')=
  \vec v'$.
  
  In this setup, it follows from the definition of $\mathcal H$ and
  Appendix~\ref{app:hompol}, Theorem~\ref{thm:ploarized_defo} that the
  section $ T\tilde \mu(s) \in H^0(\Delta \times \ell, \tilde
  \mu^*(T_X))$ is in the image of the map
  $$
  H^0(\Delta \times \ell, \tilde \mu^*(\Jet^1(L)^\vee \otimes L))
  \to H^0(\Delta \times \ell, \tilde \mu^*(T_X))
  $$
  that comes from the dualized and twisted jet
  sequence~\eqref{eq:A} of page~\pageref{eq:A}.

  To end the proof of Claim~\ref{claim:tau2}, let $\overline z' \in
  \{\pi_1(z)\} \times \ell$ be the mirror point with respect to the
  line bundle $H$. Since the coordinate change respects the line
  bundle $H$, Proposition~\ref{prop:mirrors} asserts that $\overline
  z' \in \tau'_1$. In particular, we have that $s(\overline z') \in
  T_{\tau'_1}|_{\overline z'}$ and therefore, since $\tau'_1$ is
  contracted, $T\tilde \mu (s(\overline z')) = 0$.
  Proposition~\ref{prop:section2} implies that $T\tilde \mu(s(z')) =
  0$, too. This shows that $\mu$ contracts $\tau_2$ to a point. The
  proof of Claim~\ref{claim:tau2} is finished.
\end{proof}

\begin{proof}[Application of Claim~\ref{claim:tau2}]
  Using Claim~\ref{claim:tau2}, we will derive a contradiction,
  showing that the Assumption~\ref{ass:thm11} is absurd. The proof of
  Theorem~\ref{thm:main2and3} will then be finished.
  
  For this, observe that $\tau_1 \cap \pi^{-1}(f) = \{f\} \times \{
  \sigma \}$. The sections $\tau_1$ and $\tau_2$ are therefore not
  disjoint. In this setup, Claim~\ref{claim:tau2} implies that
  $\mu(\tau_2) = \{x\}$, so that $\tau_2 \subset \mu^{-1}(x)$. That
  violates the decomposition Lemma~\ref{lem:decomp} from above.
\end{proof}

\subsection{Subbundles in the pull-back of $F$ and $T_X$}\label{sec:subb}

We will now lay the ground for the proof of
Proposition~\ref{prop:section2} in the next section. Our line of
argumentation is based on following fact which is an immediate
consequence of the Assumption~\ref{ass:thm11} and the infinitesimal
description of the universal morphism $\mu$.

\begin{fact}[{\cite[Prop.~II.3.4]{K96}, Fact~\ref{fact:tgt_vect_to_Hom}}]\label{fact:1}
  In the setup of section~\ref{sect:4setup}, let $\sigma \in \ell$ be
  the singular point. Then the restriction morphism
  $$
  H^0(\ell, f^*(T_X)) \to f^*(T_X)|_{\sigma} 
  $$
  is surjective. In other words, the vector space
  $f^*(T_X)|_{\sigma}$ is generated by global sections.
\end{fact}

Recall from Section~\ref{sect:contLines}, Fact~\ref{fact:perp} that
the non-negative part of the restriction of the vector bundle $F$ to
one of the smooth contact lines $\ell_i$ was denoted by $F|^{\geq
  0}_{\ell_i}$. We use Fact~\ref{fact:1} to show that the two vector
bundles $F|_{\ell_1}^{\geq 0}$ on $\ell_1$ and $F|_{\ell_2}^{\geq 0}$
on $\ell_2$ together give a global vector bundle on $\ell$.

\begin{lem}\label{lem:fgeq0}
  There exists a vector bundle $f^*(F)^{\geq 0} \subset f^*(F)$ on
  $\ell$ whose restriction to any of the irreducible components
  $\ell_i$ equals $F|_{\ell_i}^{\geq 0} \subset f^*(F)$. If $y \in
  \ell$ is a general point, then the restriction morphism
  $$
  H^0(\ell, f^*(F)^{\geq 0}) \to f^*(F)^{\geq 0}|_y
  $$
  is surjective.
\end{lem}
\begin{proof}
  By Fact~\ref{fact:1}, we can find sections $s_1, \ldots, s_{2n-1}
  \in H^0(\ell, f^*(T_X))$ that span
  $$
  Tf(T_{\ell_1}|_{\sigma})^\perp = Tf(T_{\ell_2}|_{\sigma})^\perp
  \subset f^*(F)|_{\sigma}
  $$
  where $\sigma \in (\ell_1 \cap \ell_2)_{\red}$ is the singular
  point of $\ell$, and where $\perp$ means: ``perpendicular with
  respect to the O'Neill tensor $N$''. Note that the sections $s_1,
  \ldots, s_{2n-1}$ become linearly dependent only at smooth points of
  the curve $\ell$.  Thus, the double dual of the sheaf generated by
  $s_1, \ldots, s_{2n-1}$ is a locally free subsheaf of $f^*(T_X)$.
  
  It follows from Fact~\ref{fact:3} that $s_1|_{\ell_i}, \ldots,
  s_{2n-1}|_{\ell_i}$ are in fact sections of $f^*(F)$ that generate
  $F|_{\ell_i}^{\geq 0}$ on an open set of $\ell_i$. The claim
  follows.
\end{proof}

\begin{cor}\label{cor:tangent}
  There exists a vector sub-bundle $T \subset f^*(F)$ whose
  restriction to any component $\ell_i$ is exactly the image of the
  tangent map
  $$
  T|_{\ell_i} = \Image(Tf :  T_{\ell_i} \to f^*(T_X)|_{\ell_i}).
  $$
  Since $f|_{\ell_i}$ is an embedding, $T|_{\ell_i}$ is of degree
  2.
\end{cor}
\begin{proof}
  By Fact~\ref{fact:perp}, we can set $T := (f^*(F)^{\geq 0})^\perp$.
\end{proof}

The vector bundle $f^*(F)^{\geq 0}$ is a sub-bundle of both $f^*(F)$
and $f^*(T_X)$. As a matter of fact, it appears as a direct summand in
these bundles.

\begin{lem}\label{lem:6}
  The vector bundle sequences on $\ell$
  \begin{equation}\label{eq:F-sequence}
    0 \to f^*(F)^{\geq 0} \to f^*(F) \to \factor f^*(F) . f^*(F)^{\geq 0}. \to 0
  \end{equation}
  and
  \begin{equation}\label{eq:TX-sequence}
    0 \to f^*(F)^{\geq 0} \to f^*(T_X) \to \factor f^*(T_X) . f^*(F)^{\geq 0}. \to 0
  \end{equation}
  are both split. We have $\factor f^*(T_X).f^*(F)^{\geq 0}. \cong
  \O_\ell \oplus \O_\ell$.
\end{lem}

\begin{proof}
  In order to show that sequence~\eqref{eq:F-sequence} splits, we show
  that the obstruction group
  $$
  \Ext^1_\ell \left(\factor f^*(F).f^*(F)^{\geq 0}., f^*(F)^{\geq 0} \right) =
  H^1 \Bigl( \ell, \underbrace{\Bigl(\factor f^*(F).f^*(F)^{\geq 0}.
    \Bigr)^\vee \otimes f^*(F)^{\geq 0}}_{=:\mathcal E} \Bigr)
  $$
  vanishes. If $\ell_i \subset \ell$ is any component, it follows
  immediately from Fact~\ref{fact:perp} that
  $$
  \left.\left( \factor f^*(F).f^*(F)^{\geq 0}. \right)\right|_{\ell_i} \cong \O_{\ell_i}(-1)
  $$
  and 
  $$
  \mathcal E|_{\ell_i} \cong \O_{\ell_i}(3) \oplus
  \O_{\ell_i}(2)^{\oplus n-1} \oplus \O_{\ell_i}(1)^{\oplus n-1}.
  $$
  By Lemma~\ref{lem:vanishingForDubbies}, $H^1(\ell, \mathcal E) =
  0$.  That shows the splitting of the sequence~\eqref{eq:F-sequence}.
  
  As a next step, we will show that the quotient $\factor
  f^*(T_X).f^*(F)^{\geq 0}.$ is trivial.  By Fact~\ref{fact:1}, we can
  find two sections $s_1, s_2 \in H^0(\ell, f^*(T_X))$ such that the
  induced sections $s'_1, s'_2 \in H^0 \left(\ell, \factor
    f^*(T_X).f^*(F)^{\geq 0}.\right)$ generate the quotient $\left.
    \factor f^*(T_X).  f^*(F)^{\geq 0}.\right|_{\sigma}$ at the
  singular point $\sigma \in \ell$.  Restricting these sections to
  $\ell_i$, it follows that the sections
  $$
  s'_1|_{\ell_i}, s'_2|_{\ell_i} \in H^0 \Bigl( \ell_i,
  \underbrace{ \left. \factor f^*(T_X).f^*(F)^{\geq
        0}.\right|_{\ell_i}}_{\cong \O_{\ell_i} \oplus \O_{\ell_i}
  \text{ by Fact~\ref{fact:perp}}}
  \Bigr)
  $$
  do not vanish anywhere and are everywhere linearly independent.
  Consequence: the induced morphism of sheaves on $\ell$
  $$
  \begin{array}{ccc}
    \O_\ell \oplus \O_\ell & \to & \factor f^*(T_X).f^*(F)^{\geq 0}. \\
    (g,h)        & \mapsto & g\cdot s'_1 + h\cdot s'_2
  \end{array}
  $$
  is an isomorphism, and the map
  $$
  \begin{array}{ccc}
    \O_\ell \oplus \O_\ell & \to     & f^*(T_X) \\
    (g,h)        & \mapsto & g\cdot s_1 + h\cdot s_2
  \end{array}
  $$
  splits the sequence~\eqref{eq:TX-sequence}.
\end{proof}

\begin{cor}\label{cor:8}
  The natural morphism 
  $$
  H^0(\ell, f^*(\Omega^1_X)) \to H^0(\ell, f^*(F)^\vee),
  $$
  which comes from the dual of the contact
  sequence~\eqref{eq:contact-sequence} of
  page~\pageref{eq:contact-sequence}, is an isomorphism.
\end{cor}
\begin{proof}
  The morphism is part of the long exact sequence 
  $$
  0 \to H^0(\ell, f^*(L)^\vee) \to H^0(\ell, f^*(\Omega^1_X)) \to
  H^0(\ell, f^*(F)^\vee)  \to \cdots
  $$
  Since $f^*(L)^\vee$ is a line bundle whose restriction to any
  irreducible component $\ell_i \subset \ell$ is of degree $-1$, there
  are no sections to it: $h^0(\ell, f^*(L)^\vee) = 0$. It remains to
  show that $h^0(\ell, f^*(\Omega^1_X))= h^0(\ell, f^*(F)^\vee)$. The
  direct sum decomposition of Lemma~\ref{lem:6} yields
  $$
  \begin{array}{rcl}
    h^0(\ell, f^*(F)^\vee) & = & h^0(\ell, (f^*(F)^{\geq 0})^\vee) +
    \underbrace{h^0 \Bigl(\ell, \Bigl( \factor f^*(F).f^*(F)^{\geq
    0}. \Bigr)^\vee \Bigr)}_{= 2 \text{ by Fact~\ref{fact:perp} and Proposition~\ref{prop:mirrors}}} \\
    h^0(\ell, f^*(\Omega^1_X)) & = & h^0(\ell, (f^*(F)^{\geq 0})^\vee) +
    \underbrace{h^0 \Bigl(\ell, \Bigl(\factor f^*(T_X).f^*(F)^{\geq
    0}. \Bigr)^\vee \Bigr)}_{= h^0(\O_\ell \oplus \O_\ell) = 2 \text{
    by Lemma~\ref{lem:6}}} 
  \end{array}
  $$
  The corollary follows.
\end{proof}

\subsection{The vanishing locus of sections in the pull-back of $T_X$}\label{sec:secinsubb}

Using Corollary~\ref{cor:8}, we can now establish a criterion,
Proposition~\ref{prop:section2}, that guarantees that certain sections
in $f^*(T_X)$ that vanishes at a point $y \in \ell$ will also vanish
at the mirror point. The following lemma is a first precursor.

\begin{lem}\label{lem:section1}
  In the setup of section~\ref{sect:4setup}, let $y \in \ell$ be a
  general point and let $s \in H^0(\ell, f^*(T_X))$ be a section that
  vanishes at $y$. Then the associated section $s' \in H^0 \left(\ell,
    \factor f^*(T_X). T. \right)$ vanishes at the mirror point
  $\overline y$.  Here $T$ is the vector bundle from
  Corollary~\ref{cor:tangent}.
\end{lem}
\begin{proof}
  We claim that $s \in H^0(\ell, f^*(F))$. The proof of this claim is
  a twofold application of Fact~\ref{fact:3}. If we assume without
  loss of generality that $y \in \ell_1$, then a direct application of
  Fact~\ref{fact:3} shows that $s|_{\ell_1} \in H^0(\ell_1,
  F|_{\ell_1}^{\geq 0})$. If $\sigma = (\ell_1 \cap \ell_2)_{\red}$ is
  the singular point, this implies that $s(\sigma) \in
  (F|_{\ell_1}^{\geq 0})|_\sigma = (F|_{\ell_2}^{\geq 0})|_\sigma$.
  Another application of Fact~\ref{fact:3} then shows the claim.
 
  Consequence: in order to show Lemma~\ref{lem:section1} it suffices
  to show that the associated section $s'' \in H^0 \left(\ell, \factor
    f^*(F).T.  \right)$ vanishes at $\overline y$. We assume to the
  contrary.
  
  Since $T^\perp = f^*(F)^{\geq 0}$, the non-degenerate O'Neill tensor
  gives an identification
  $$
  f^*(F)^{\geq 0}|_y \cong \left. \Bigl( \Bigl (\factor f^*(F).T.
    \Bigr)^\vee \otimes f^*(L)\Bigr)\right|_y
  $$
  By Lemma~\ref{lem:fgeq0}, we can therefore find a section $t \in
  H^0(\ell, f^*(F)^{\geq 0})$ such that $s$ and $t$ pair to give a
  section
  $$
  N(s, t) \in H^0(\ell, f^*(L))
  $$
  That vanishes at $y$, but does not vanish on $\overline y$.  That
  is a contradiction to Proposition~\ref{prop:mirrors}.
\end{proof}

In Lemma~\ref{lem:section1}, it is generally \emph{not} true that the
section $s$ vanishes at $\overline y$ ---to a given section $s$, we
can always add a vector field on $\ell$ that stabilizes $y$, but does
not stabilize the mirror point $\overline y$. However, the statement
becomes true if we restrict ourselves to sections $s$ that come from
$L$-jets.

\begin{prop}\label{prop:section2}
  In the setup of section~\ref{sect:4setup}, let $y \in \ell$ be a
  general point and $s \in H^0(\ell, f^*(T_X))$ a section that
  vanishes on $y$. If $s$ is in the image of the map
  $$
  H^0(\ell, f^*(\Jet^1(L)^\vee \otimes L)) \to H^0(\ell, f^*(T_X)),
  $$
  that comes from the dualized and twisted jet
  sequence~\eqref{eq:B}, then $s$ vanishes also at the mirror point
  $\overline y$.
\end{prop}

The proof of Proposition~\ref{prop:section2} requires the following
Lemma, which we state and prove first.

\begin{lem}\label{lem:JetCriterion2}
  Let $s \in H^0(\ell, f^*(F))$ be a section and let $D \in |f^*(L)|$
  be an effective divisor that is supported on the smooth locus of
  $\ell$. If $s$ vanishes on $D$, then $s$ is in the image of the map
  $$
  H^0(\ell, f^*(\Jet^1(L)^\vee\otimes L)) \to H^0(\ell, f^*(T_X))
  $$
  that comes from the dualized and twisted jet
  sequence~\eqref{eq:B} of page~\pageref{eq:B}.
\end{lem}
\begin{proof}
  In view of Fact~\ref{fact:jetsoncontact}, we need to show that
  $s$ is in the image of the map
  $$
  \alpha : H^0(\ell, f^*(\Omega^1_X\otimes L)) \to H^0(\ell,
  f^*(F))
  $$
  which comes from the dualized and twisted contact
  sequence~\eqref{eq:A}. For that, let $t \in H^0(\ell, f^*(L))$ be a
  non-zero section that vanishes on $D$. Using the O'Neill tensor $N$
  to identify $F$ with $F^\vee \otimes L$, we can view $s$ as a
  section that lies in the image
  $$
  \begin{CD}
    H^0(\ell, f^*(F^\vee)) @>{\cdot t}>> H^0(\ell, f^*(F^\vee \otimes L))
  \end{CD}
  $$
  The claim then follows from Corollary~\ref{cor:8}, and the
  commutativity of the diagram
  $$
  \xymatrix{
    {H^0(\ell, f^*(\Omega^1_X))} \ar[r]^{\text{surjective}}_{\text{by
        Cor.~\ref{cor:8} }} \ar[d]_{\cdot t}&
    H^0(\ell, f^*(f^\vee))  \ar[d]^{\cdot t} \\
    {H^0(\ell, f^*(\Omega^1_X \otimes L))} \ar[r] & {H^0(\ell,
      f^*(F^\vee \otimes L))} }
  $$
\end{proof}

\begin{proof}[Proof of Proposition~\ref{prop:section2}]
  Since $s \in H^0(\ell, f^*(F))$, Fact~\ref{fact:jetsoncontact}
  implies that $s$ is in the image of the map $\alpha$ from the long
  exact sequence associated with the dualized and twisted
  Contact-sequence~\eqref{eq:A}
  \begin{multline}\label{eq:longExRedJet}
    0 \to \underbrace{H^0(\ell, f^*(\O_X))}_{\cong \mathbb C} \to
    H^0(\ell, f^*(\Omega^1_X\otimes L)) \xrightarrow{\alpha} H^0(\ell,
    f^*(F)) \to \\
    \to \underbrace{H^1(\ell, f^*(\O_X))}_{\cong \mathbb C \text{ by
        Remark~\ref{rem:dubbyconic}}} \to \cdots
  \end{multline}
  By Lemma~\ref{lem:section1}, the vector space
  $$
  H_{y\overline y} := \{ \tau \in H^0(\ell, f^*(F))\,|\, \tau(y)=0, \,
  \tau(\overline y)=0\}
  $$
  is a linear hyperplane in
  $$
  H_y := \{ \tau \in H^0(\ell, f^*(F))\,|\, \tau(y)=0 \}.
  $$
  Because $\O_\ell(y+\overline y) \cong f^*(L)$,
  Lemma~\ref{lem:JetCriterion2} implies that
  $$
  H_{y \overline y} \subset \Image(\alpha) \cap H_y.
  $$
  But $\codim_{H^0(\ell, f^*(F))} \Image(\alpha) \leq 1$, so that
  there are only two possibilities here:
  \begin{enumerate}
  \item $H_y \subseteq \Image(\alpha)$ and $\Image(\alpha) \cap H_y =
    H_y$
  \item $\Image(\alpha) \cap H_y  = H_{y\overline y}$
  \end{enumerate}
  Observe that Proposition~\ref{prop:section2} is shown if we rule out
  possibility (1). For that, it suffices to show that there exists a
  section $t \in H_y$ which is not in the image of $\alpha$.
  
  To this end, let $\theta \in H^0(X, \Omega^1_X \otimes L)$ be the
  nowhere-vanishing $L$-valued 1-form that defines the contact
  structure in Sequence~\eqref{eq:contact-sequence} of
  page~\pageref{eq:contact-sequence}. The beginning part of
  Sequence~\eqref{eq:longExRedJet} says that its pull-back
  $f^*(\theta) \in H^0(\ell, f^*(\Omega^1_X \otimes L))$ is, up to
  multiple, the unique section that is in the kernel of $\alpha$. If
  we fix $i \in \{0,1\}$, then the analogous sequence for
  $f|_{\ell_i}$ tells us that $f^*(\theta)|_{\ell_i}$ is the unique
  (again up to a multiple) section in $H^0(\ell_i,
  (f|_{\ell_i})^*(\Omega^1_X \otimes L))$ which is in the kernel of
  $\alpha|_{\ell_i}$. Consequence: there exists no section $u \in
  H^0(\ell, f^*(\Omega^1_X \otimes L))$ such that $\alpha(u)$ vanishes
  on one component of $\ell = \ell_1 \cup \ell_2$, but not on the
  other.

  By Lemma~\ref{lem:bideg22}, however, there exists a section $t \in
  H^0(\ell, T) \subset H^0(\ell, f^*(F))$ that vanishes on the
  component of $y$ and not on the other. The section $t$ is therefore
  contained in $H_y$ but not in $\Image(\alpha)$. This ends the proof
  of Proposition~\ref{prop:section2}.
\end{proof}

\section{Contact lines sharing more than one point}\label{chapt:banana}

As a last step before the proof of the main theorem, we show property
\iref{l4} from the list of Theorem~\ref{thm:main}.

\begin{thm}\label{thm:banana}
  Let $x \in X$ be a general point and let $\ell_1$, $\ell_2$ be two
  distinct contact lines through $x$. Then $\ell_1$ and $\ell_2$
  intersect in $x$ only, $\ell_1 \cap \ell_2 = \{x\}$.
\end{thm}

The proof is really a corollary to the results of the previous
chapter. In analogy to Definition~\ref{def:dubby}, we name the
simplest arrangement of rational curves that intersect in two points.

\begin{defn}
  A \emph{pair with proper double intersection} is a reduced,
  reducible curve, isomorphic to the union of a line and a smooth
  conic in $\P^2$ intersecting transversally in two points.
\end{defn}

\begin{proof}[Proof of Theorem~\ref{thm:banana}]
  We argue by contradiction and assume that for a general point $x$
  there is a pair of contact lines $\ell_1$, $\ell_2$ through $x$
  which meet in at least one further point. The pair $\ell_1 \cup
  \ell_2$ will then be dominated by a pair with proper double
  intersection $\ell = \ell'_1 \cup \ell'_2$. More precisely, there
  exists a generically injective morphism $f: \ell \to \ell_1 \cup
  \ell_2$ which maps $\ell'_i$ to $\ell_i$ and which maps one of the
  two singular points of $\ell$ to $x$.  Let $y \in \ell$ be that
  point.
  
  Since $x$ is assumed to be a general point, there exists an
  irreducible component of the reduced $\Hom$-scheme
  $$
  \mathcal H \subset [\Hom(\ell, X)]_{\red}
  $$
  with universal morphism $\mu: \mathcal H \times \ell \to X$ such
  that the restriction
  $$
  \mu' = \mu|_{\mathcal H \times \{y\}} : \mathcal H \to X
  $$
  is dominant. We can further assume that $f$ is a smooth point of
  $\mathcal H$ and that the tangent map $T\mu'$ has maximal rank
  $2n+1$ at $f$.
  
  By Fact~\ref{fact:perp} and Theorem~\ref{thm:main2and3}, the tangent
  spaces $T_{\ell_1}|_x \subset F|_x$ and $T_{\ell_2}|_x \subset
  F|_x$ are both 1-dimensional and distinct. We can thus find a
  tangent vector $\vec v \in F|_x$ which is perpendicular (with
  respect to the non-degenerate O'Neill-tensor) to $T_{\ell_1}|_x$ but
  not to $T_{\ell_2}|_x$. Since the tangent map
  $$
  T \mu'|_f : T_{\mathcal H}|_f \to T_X|_x
  $$
  has maximal rank, we can find a tangent vector $s \in T_{\mathcal
    H}|_f$ such that $T\mu'(s) = \vec v$. By \cite[II.~Prop.~3.4]{K96}
  that means that we can find a section
  $$
  s \in T_{\mathcal H}|_f = H^0(\ell, f^*(T_X))
  $$
  with $T\mu'(s(y)) = \vec v \in f^*(T_X)$.
  
  Now let $\theta: T_X \to L$ be the $L$-valued 1-form that defines
  the contact structure in Sequence~\eqref{eq:contact-sequence} of
  page \pageref{eq:contact-sequence}. We need to consider the section
  $s' := f^*(\theta)(s) \in H^0(\ell, f^*(L))$. Recall
  Fact~\ref{fact:3} which asserts that $s'$ vanishes identically on
  $\ell'_1$, but does not vanish identically on $\ell'_2$. In
  particular, $\ell'_1 \cap \ell'_2$ is contained in the zero-locus of
  $s'|_{\ell'_2}$ and we have
  $$
  \deg f^*(L)|_{\ell'_2} \geq \#(\ell'_1 \cap \ell'_2) = 2.
  $$
  But $\ell_2$ is a contact line and $f^*(L)|_{\ell'_2}$ is a line
  bundle of degree 1, a contradiction.
\end{proof}

\section{Proof of the main results}\label{chapt:mainthmproof}

\subsection{Proof of Theorem~\ref{thm:main}}
In view of Theorem~\ref{thm:irreducibility}, to prove
Theorem~\ref{thm:main}, it only remains to show that $\locus(H_x)$ is
a cone. This will turn out to be a corollary to
Theorems~\ref{thm:main2and3} and \ref{thm:banana}.

Let $\tilde H_x$ be the normalization of the subspace $H_x \subset H$
of contact line through $x$. Since all contact lines through $x$ are
free, it follows from \cite[Chapt.~II, Prop.~3.10 and
Cor.~3.11.5]{K96} that $\tilde H_x$ is smooth. We have a diagram
$$
\xymatrix{
  &&& {\hat X} \ar[d]_{\beta}^{\text{blow-up of $x$}} \\
  {\tilde U_x} \ar[rrr]^{\iota}_{\text{evaluation morphism}}
  \ar[d]_{\pi} \ar@/^/@{-->}[rrru]^{\hat \iota = \beta^{-1} \circ
    \iota}
  &&& X \\
  {\tilde H_x} }
$$
where $\tilde U_x$ is the pull-back of the universal $\P^1$-bundle
$\Univrc(X)$, $\iota$ the natural evaluation morphism, and $\hat X =
\BlowUp(X,x)$ the blow-up of $x$ with exceptional divisor $E$. Since
all contact lines through $x$ are smooth, the scheme-theoretic fiber
$\iota^{-1}(x)$ is a Cartier-divisor in $\tilde U_x$, and it follows
from the universal property \cite[Chapt. II, Prop.~7.14]{Ha77} of the
blow-up that $\hat \iota = \beta^{-1} \circ \iota$ is actually a
morphism.
  
To show that $\locus(H_x) = \Image(\iota)$ really is a cone in the
sense of \cite[Chapt.~1.1.8]{BS95}, it suffices to show that $\hat
\iota$ is an embedding, i.e., that $\hat \iota$ is injective and
immersive.

\begin{description}
\item[Injective] Let $y \in \Image(\hat \iota)$ be any point. If $y
  \in E$, Theorem~\ref{thm:main2and3} asserts that $\# \hat
  \iota^{-1}(y)=1$. If $y \not \in E$, the same is guaranteed by
  Theorem~\ref{thm:banana}.
    
\item[Immersive] Fact~\ref{fact:perp} implies that for every
  $\pi$-fiber $\ell \cong \P^1$, we have
  $$
  \hat \iota^*(T_{\hat X})|_{\ell} \cong \O_{\P^1}(2) \oplus
  \O_{\P^1}^{\oplus n-1} \oplus \O_{\P^1}(-1)^{\oplus n+1}.
  $$
  Under this condition, \cite[Chapt. II, Prop.~3.4]{K96} shows that
  $\hat \iota$ is immersive as required.
\end{description}
This ends the proof of Theorem~\ref{thm:main}. \qed

\subsection{Proof of Corollary~\ref{cor:stability}}

Once Theorem~\ref{thm:main} is shown, the proof of
\cite[Thm.~2.11]{Hwa00} applies nearly verbatim to contact manifolds
---note, however, that the estimate of \cite[Thm.~2.11]{Hwa00} is not
optimal.  For the reader's convenience, we recall the argumentation
here.

Assume that the tangent bundle $T_X$ is not stable. By
\cite[Prop.~4]{Hwa98}, this implies that we can find a subsheaf
$\mathcal G \subset T_X$ of positive rank with the following
intersection property. If $x \in X$ is a general point, $\mathcal C_x
\subset \P(T_X|^\vee_x)$ the projective tangent cone of $\locus(H_x)$,
$y \in \mathcal C_x$ a general point and $T \subset \P(T_X|^\vee_x)$
the projective tangent space to $\mathcal C_x$ at $y$,
then
\begin{equation}  \label{eq:intersect}
\dim (T \cap \P(\mathcal G|^\vee_x)) \geq \frac{\rank (\mathcal
  G)}{\dim X}(n+1)-1.
\end{equation}
We will show that this leads to a contradiction. Let
$$
\psi: \P(T_X|^\vee_x) \setminus \P(\mathcal G|^\vee_x) \to \P^{\dim X
  -\rank (\mathcal G)-1}
$$
be the projection from $\P(\mathcal G|^\vee_x)$ to a complementary
linear space, and let $q$ be the generic fiber dimension of
$\psi|_{\mathcal C_x}$. We will give two estimates for $q$.

\subsubsection*{Estimate 1} 
Since a tangent vector in $T_{\mathcal C_x}|_y$ is in the kernel of
the tangent map $T(\psi|_{\mathcal C_x})$ if the associated line in
$T$ intersects $\P(\mathcal G|^\vee_x)$, equation~\eqref{eq:intersect}
implies that the kernel of $T(\psi|_{\mathcal C_x})$ is of dimension
$$
  \dim \ker( T(\psi|_{\mathcal C_x})) \geq \frac{\rank (\mathcal G)}{\dim X}(n+1)
$$
Consequence:
\begin{equation}\label{eq:qestimate}
  q \geq \frac{\rank (\mathcal G)}{\dim X}(n+1).
\end{equation}

\subsubsection*{Estimate 2}
Let $T' \subset \P^{\dim X - \rank (\mathcal G) - 1}$ be the
projective tangent space to the (smooth) point $\psi(y)$ of the image
of $\psi$. Then $\psi^{-1}(T')$ is a linear projective subspace of
dimension
$$
\dim \psi^{-1}(T') = \dim T +\rank (\mathcal G) = 
(\dim \mathcal C_x -q)+\rank (\mathcal G).
$$
This linear space is tangent to $\mathcal C_x$ along the fiber of
$\psi|_{\mathcal C_x}$ through $y$. Since $\mathcal C_x$ is smooth by
Theorem~\ref{thm:main}, Zak's theorem on tangencies, \cite{Zak93} (see
also \cite[Thm.~2.7]{Hwa00}), asserts that
\begin{align*}
  & & \dim (\text{fiber of $\psi|_{\mathcal C_x}$ through $\alpha$})
  & \leq \dim (\psi^{-1}(T')) - \dim \mathcal C_x \\
  \Rightarrow && q & \leq (\dim \mathcal C_x-q+\rank (\mathcal G)) - \dim \mathcal C_x \\
  \Rightarrow && q & \leq \frac{\rank (\mathcal G)}{2}
\end{align*}

\subsubsection*{Application of the Estimates}
Combining Estimate~2 with \eqref{eq:qestimate}, we obtain
\begin{align*}
  && \frac{\rank (\mathcal G)}{\dim X}(n+1) & 
  \leq \frac{\rank (\mathcal G)}{2} \\
  \Rightarrow && 2(n+1) & \leq \dim X
\end{align*}
But we have $\dim X = 2n+1$, a
contradiction. Corollary~\ref{cor:stability} is thus shown. \qed

\subsection{Proof of Corollary~\ref{cor:extension}}

This corollary follows directly from Theorem~\ref{thm:main} and
\cite[Thm.~3.2]{Hwa00}.

\appendix

\section{A description of the jet sequence}
\label{app:jets}

\subsection{The first jet sequence}
\label{app:jetsequence}
Let $X$ be a complex manifold (not necessarily projective or compact)
and $L \in \Pic(X)$ a line bundle.  Throughout the present paper, we
use the definition for the first jet bundle $\Jet^1(L)$ that was
introduced by Kumpera and Spencer in \cite[Chapt.~2]{KS72} and now
seems to be standard in algebraic geometry ---see also
\cite[Chapt.~1.6.3]{BS95}.  One basic feature of the first jet bundle
of $L$ is the existence of a certain sequence of vector bundles, the
first jet sequence of $L$.
\begin{equation}\label{eq:firstJetSeq}
  0 \longrightarrow \Omega^1_X \otimes L \xrightarrow{\gamma} \Jet^1(L)
  \xrightarrow{\delta} L \longrightarrow 0
\end{equation}
There exists a morphism of sheaves,
$$
\Prolong : L \to \Jet^1(L),
$$
called the ``prolongation'' which makes \eqref{eq:firstJetSeq} a
split sequence of sheaves. The first jet sequence is, however,
generally \emph{not} split as a sequence of vector bundles, and the
prolongation morphism is definitely not $\O_X$-linear. In fact, an
elementary computation using the definition of $\Jet^1(L)$ from
\cite{KS72} and the construction of differentials from
\cite[Chapt.~25]{Mat89} yields that for any open set $U \subset X$,
any section $\sigma \in L(U)$ and function $g \in \O_X(U)$, we have
$$
\Prolong(g\cdot \sigma) = g \cdot \Prolong(\sigma) + \gamma ( dg
\otimes \sigma).
$$

\subsection{Jets and logarithmic differentials}

The definitions of \cite{KS72} are well suited for algebraic
computations. If we are to apply jets to deformation-theoretic
problems, however, it seems more appropriate to follow an approach
similar to that of Atiyah, \cite{Ati57}, and to describe jets in terms
of logarithmic tangents and differentials on the (projectivized) total
space of the line bundle.  We refer to \cite[Chapt.~2.1]{KPSW00} for a
brief review of Atiyah's definitions.  While the relation between
\cite{KS72} and our construction here is probably understood by
experts, the author could not find any reference. A detailed
description is therefore included here.

Set $Y := \P(L \oplus \O_X)$. We denote the natural $\P^1$-bundle
structure by $\pi : Y \to X$ and let $\Sigma = \Sigma_0 \cup
\Sigma_\infty \subset Y$ be the union of the two disjoint sections
that correspond to the direct sum decomposition. By convention, let
$\Sigma_\infty$ the section whose complement $Y \setminus
\Sigma_\infty$ is canonically isomorphic to the total space of the
line bundle $L$.

Let $\Omega^1_Y(\log \Sigma)$ be the locally free subsheaf of
differentials with logarithmic poles along $\Sigma$. This sheaf, which
contains $\Omega^1_Y$ as a subsheaf, is defined and thoroughly
discussed in \cite[Chapt.~II.3]{Deligne70}. In particular, it is shown
in \cite[Chapt.~II.3.3]{Deligne70} that the sequence of relative
differentials
$$
0 \longrightarrow \pi^* \Omega^1_X  \longrightarrow \Omega^1_Y
\longrightarrow \Omega^1_{Y|X} \longrightarrow 0
$$
restricts to an exact vector bundle sequence of logarithmic
tangents
$$
0 \longrightarrow \pi^* \Omega^1_X \longrightarrow \Omega^1_Y(\log \Sigma)
\longrightarrow \Omega^1_{Y|X}(\log \Sigma) \longrightarrow 0.
$$
Since $R^1 \pi_* (\pi^*(\Omega^1_X)) = 0$, we can push down to $X$,
twist by $L$ and obtain a short exact sequence as follows
\begin{equation}\label{eq:rellogs}
  0 \longrightarrow \Omega^1_X \otimes L \xrightarrow{\beta} \pi_*
  \Omega^1_Y(\log \Sigma) \otimes L
  \longrightarrow \pi_* \Omega^1_{Y|X}(\log \Sigma) \otimes L \longrightarrow 0.
\end{equation}
We will show that sequence~\eqref{eq:rellogs} is canonically
isomorphic to the first jet sequence~\eqref{eq:firstJetSeq} of $L$.

\begin{thm}\label{thm:jetsandlogs1}
  With the notation from above, there exists an isomorphism of vector bundles 
  $$
  \alpha : \pi_* \Omega^1_Y(\log \Sigma) \otimes L \to \Jet^1(L)
  $$
  such that the diagram
  \begin{equation}\label{diag:JetsAndLogs}
    \xymatrix{
      0 \ar[r] & {\Omega^1_X\otimes L} \ar[r]^(0.35){\beta} 
      \ar[d]_{\text{Identity}} & {\pi_* \Omega_Y^1(\log \Sigma)
      \otimes L} \ar[r] \ar[d]^{\alpha} &
      {\pi_* \Omega^1_{Y|X}(\log \Sigma) \otimes L} \ar[r] & 0 \\
      0 \ar[r] & {\Omega^1_X\otimes L} \ar[r]^{\gamma} & {\Jet^1(L)}
      \ar[r]^{\delta} & {L} \ar[r] \ar@/^0.3cm/[l]^{\Prolong} & 0
    }
  \end{equation}
  commutes, i.e.~$\gamma = \alpha \circ \beta$.
\end{thm}

In the following Appendix~\ref{app:hompol}, where deformations of
morphisms are discussed, we will need to consider tangents rather than
differentials. For that reason, we state a ``dualized and twisted''
version of Theorem~\ref{thm:jetsandlogs1}.  Recall from
\cite[Chapt.~II.3]{Deligne70} that the dual of $\Omega^1_Y(\log
\Sigma)$ is the locally free sheaf $T_Y(-\log \Sigma)$ of vector
fields on $Y$ which are tangent to $\Sigma$.

\begin{cor}\label{cor:jetsandlogs}
  There exists an isomorphism $A$ of vector bundles such that the
  diagram
  $$
  \xymatrix{
    {\Jet^1(L)^\vee \otimes L} \ar[rr]^(.6){\text{tangent map $T\pi$}} 
    \ar[d]^{A} & & {T_X} \ar[d]^{\text{Identity}} \\
    {\pi_* T_Y(-\log \Sigma)} \ar[rr]_(.6){\gamma^\vee \otimes Id_L}
    & & {T_X}
  }
  $$
  commutes. \qed
\end{cor}

Informally speaking, we can say the following.

\begin{summary}
  A vector field on the manifold $X$ comes from $L$-jets if and only
  if it lifts to a vector field on $Y$ whose flow stabilizes
  $\Sigma_0$ and $\Sigma_\infty$.
\end{summary}

\subsubsection*{Proof of Theorem~\ref{thm:jetsandlogs1}, setup}

Let $U \subset X$ be an open set and $\sigma \in L(U)$ a
nowhere-vanishing section. We will construct the isomorphism $\alpha$
locally at first by defining an $\O_X$-linear morphism
$$
\alpha_{U,\sigma} : [\pi_* \Omega^1_Y(\log \Sigma) \otimes L](U) \to
[\Jet^1(L)](U)
$$
which we will later show to not depend on the choice of the section
$\sigma$. It will then follow trivially from the construction that the
various $\alpha_{U,\sigma}$ glue together to give a morphism of vector
bundles.

Throughout the proof of Theorem~\ref{thm:jetsandlogs1}, we constantly
identify sections $[\pi_* \Omega^1_Y(\log \Sigma) \otimes L](U)$ with
$[\Omega^1_Y(\log \Sigma) \otimes \pi^*(L)](\pi^{-1}(U))$. Likewise,
we will use the letter $\sigma$ to denote the subvariety of
$\P(L\oplus \O_X)|_U$ that is associated with the section.

\subsubsection*{Proof of Theorem~\ref{thm:jetsandlogs1}, definition of $\alpha_{U,\sigma}$} 

In order to define $\alpha_{U,\sigma}$, use the nowhere-vanishing
section $\sigma$ to introduce a bundle coordinate on $\pi^{-1}(U)$,
which we can view as a meromorphic function $z$ on $\pi^{-1}(U)$ with
a single zero along $\Sigma_0$ and a single pole along $\Sigma_\infty$
such that
$$
\pi \times z: \pi^{-1}(U) \to U\times \P^1
$$
is an isomorphism with $z|_{\sigma} \equiv 1$. The coordinate $z$
immediately gives a differential form
$$
d\log z := \frac{1}{z}dz \in [\Omega^1_Y(\log \Sigma)](\pi^{-1}(U))
$$
with logarithmic poles along both components of $\Sigma$. Note that
$d\log z$ yields a nowhere-vanishing section of the line bundle
$\Omega^1_{Y|X}(\log \Sigma)$ of relative logarithmic differentials.
Consequence: there exists a relative vector field
$$
\vec v_z \in [ T_{Y|X}(-\log \Sigma)](\pi^{-1}(U))
$$
with zeros along $\Sigma$ which is dual to $d\log z$, i.e., $(d\log
z)(\vec v_z) =1$. In the literature, $\vec v_z$ is sometimes denoted
by $z \frac{\partial}{\partial z}$, but we will not use this notation
here.

With these notations, if $\omega \in [\Omega^1_Y(\log \Sigma) \otimes
\pi^*(L)](\pi^{-1}(U))$ is a $\pi^*(L)$-valued logarithmic form, set
$$
\alpha_{U,\sigma}(\omega) := \gamma ( \underbrace{ \omega - d \log
  z \otimes \omega ( \vec v_z )}_{=: \theta} ) + z \circ \omega ( \vec
v_z ) \cdot \Prolong(\sigma).
$$
Explanation: we point out that $\omega ( \vec v_z)$ is a section of
$[\pi^*(L)](\pi^{-1}(U))$ so that we can regard $z \circ \omega ( \vec
v_z )$ as a function. It is an elementary calculation in coordinates
to see that $\theta$ is a regular $L$-valued 1-form on $\pi^{-1}(U)$
that vanishes on relative tangents. We can therefore see $\theta$ as
the pull-back of a uniquely determined $L$-valued 1-form on $U$. In
particular, $\gamma(\theta)$ is a well-defined 1-jet in
$[\Jet^1(L)](U)$.

\subsubsection*{Proof of Theorem~\ref{thm:jetsandlogs1}, injectivity}

It follows immediately from the definition that $\alpha_{U,\sigma}$ is
injective. Namely, if $\alpha_{U,\sigma}(\omega) = 0$, then the
exactness of the second row of diagram~\eqref{diag:JetsAndLogs}
implies that 
$$
\theta = 0 \text{, i.e.~that } \omega = \omega(\vec v_z)d\log (z)
$$
and
$$
z \circ \omega(\vec v_z) = 0\text{, i.e.~that } \omega(\vec v_z) = 0
$$
Together this implies that $\omega = 0$.

\subsubsection*{Proof of Theorem~\ref{thm:jetsandlogs1}, coordinate change}

Let $\tau \in L(U)$ be another nowhere-vanishing section, $\tau =
g\cdot \sigma$ with $g \in \O^*_X(U)$. The section $\tau$ gives rise
to a new bundle coordinate $z'$. We have
$$
z' = \frac{1}{g}\cdot z, \qquad
d \log z' = d\log z - d\log g
$$
and therefore
$$
\vec v_{z'} = \vec v_z.
$$
Using these equalities, it is a short computation to see that
$\alpha_{U,\sigma}$ and $\alpha_{U,\tau}$ agree:
\begin{equation*}
  \begin{split}
    \alpha_{U,\tau}(\omega) & = \gamma ( \omega - d\log z' \otimes
    \omega (\vec v_{z'} ) ) -
    z' \circ \omega ( \vec v_{z'} ) \cdot \Prolong(\tau) \\
    & = \gamma ( \omega - [ d\log z - d\log g] \otimes \omega (\vec
    v_z) ) - \\
    & \qquad \qquad \qquad \qquad \qquad  \frac{z}{g} \circ \omega (
    \vec v_z ) \cdot
    [ g\cdot \Prolong(\sigma) + \gamma(dg \otimes \sigma)] \\
    & = \gamma ( \omega - d\log z \otimes \omega( \vec v_z ) ) -
    z \circ \omega ( \vec v_z) \cdot \Prolong(\sigma) + \\
    & \qquad \qquad \qquad \qquad \qquad \gamma( d\log g \otimes
    \omega( \vec v_z ) ) -
    \frac{z}{g} \circ \omega(\vec v_z) \gamma( dg \otimes \sigma)\\
    & = \alpha_{U,\sigma}(\omega) + \gamma \biggl( \biggl[
    \underbrace{d\log g - \frac{1}{g}dg}_{=0} \biggr] \otimes
    \omega(\vec v_z) \biggr)
  \end{split}
\end{equation*}
We have thus constructed an injective morphism of sheaves. We will
later see that $\alpha$ is an isomorphism.

\subsubsection*{Proof of Theorem~\ref{thm:jetsandlogs1}, commutativity
  of Diagram~\eqref{diag:JetsAndLogs}}

Let $\theta \in [\Omega^1_X\otimes L](U)$. The image $\beta(\theta)$
is nothing but the pull-back of $\theta$ to $\pi^{-1}(U)$. In
particular, if $z$ is any bundle coordinate, we have that
$\beta(\theta)(\vec v_z) \equiv 0$. Therefore
\begin{equation*}
  \begin{split}
    \alpha_{U,\sigma}\circ \beta(\theta) & = \gamma \biggl(
    \beta(\theta) - d\log z \otimes \underbrace{\beta(\theta) ( \vec
      v_z )}_{= 0} \biggr) - z \circ \underbrace{\beta(\theta) ( \vec
      v_z )}_{=0} \cdot
    \Prolong(\sigma) \\
    & = \gamma(\beta(\theta))
  \end{split}
\end{equation*}
where we again identify a form $\theta$ with its pull-back.

\subsubsection*{Proof of Theorem~\ref{thm:jetsandlogs1}, end of proof}

It remains to show that the sheaf-morphism $\alpha$ is isomorphic,
i.e.~surjective.  Because Diagram~\eqref{diag:JetsAndLogs} is
commutitative, to show that $\alpha$ is surjective, it suffices that
$\delta \circ \alpha$ is surjective. Let $\sigma \in L(U)$ again be a
nowhere-vanishing section and let $\tau \in L(U)$ be any section,
$\tau = g\cdot \sigma$, where $g \in \O_X(U)$. We show that $\tau$ is
in the image of $\delta \circ \alpha_{U,\sigma}$.

For this, let $z$ be the bundle coordinate that is associated with
$\sigma$ and set
$$
\omega := d\log z \otimes (g \cdot \sigma )
$$
We have
\begin{equation*}
  \begin{split}
    \delta \circ \alpha_{U,\sigma}(\omega) & = \delta \bigl(
    \underbrace{\gamma(\cdots)}_{=0}+ z\circ \omega( \vec v_z ) \cdot
    \Prolong(\sigma) \bigr) \\
    & = \delta \bigl( \underbrace{z \circ \sigma}_{\equiv 1}\cdot g
    \cdot \underbrace{d \log z (\vec v_z )}_{=1} \cdot
    \Prolong(\sigma) \bigr) \\
    &= g\cdot \delta(\Prolong(\sigma)) = g\cdot \sigma = \tau.
  \end{split}
\end{equation*}
The proof of Theorem~\ref{thm:jetsandlogs1} is thus finished.\qed

\section{Morphisms between polarized varieties}
\label{app:hompol}

\subsection{The tangent space to the Hom-scheme}

Let $X$ be a complex projective manifold, $\ell$ a projective variety
and $f: \ell \to X$ a morphism. It is well-known that there exists a
scheme $\Hom(\ell,X)$ that represents morphisms $\ell \to X$ ---see
e.g.~\cite[Chapt.~I]{K96}. In particular, there exists a functorial
1:1-correspondence between closed points of $\Hom(\ell,X)$ and actual
morphisms. As a consequence we have a ``universal morphism''
$\Hom(\ell,X)\times \ell \to X$. It is known that the tangent space to
$\Hom(\ell,X)$ is naturally identified with the space of sections in
the pull-back of the tangent bundle
$$
T_{\Hom(\ell,X)}|_f \cong H^0(\ell, f^*(T_X)).
$$
In the most intuitive setup, this identification takes the following
form:
\begin{fact}\label{fact:tgt_vect_to_Hom}
  Let $\Delta$ be the unit disc with coordinate $t$ and let
  $$
  \begin{array}{rccc}
    f : & \Delta & \to & \Hom(\ell,X) \\
    & t & \mapsto & f_t
  \end{array}
  $$
  be a family of morphisms. If $\mu: \Delta\times \ell \to X$ is the
  induced universal morphism, then 
  $$
  Tf \left( \left. \frac{\partial}{\partial t}\right|_{t=0}\right) \in
  f^*(T_{\Hom(\ell,X)}|_{f_0}) \cong T_{\Hom(\ell,X)}|_{f_0}
  $$
  is naturally identified with
  $$
  \left. T\mu \left( \frac{\partial}{\partial t}\right)
  \right|_{\{0\}\times \ell} \in H^0(\{0\}\times \ell,
  \mu^*(T_X)|_{\{0\}\times \ell}) \cong H^0( \ell, f_0^*(T_X))
  $$
  \qed
\end{fact}

The identification has become so standard that we often wrongly write
``equal'' rather than ``naturally isomorphic''.

\subsection{The pull-back of line bundles}

In this paper we need to consider morphisms of polarized varieties.
More precisely, we fix a line bundle $L \in \Pic(X)$ and wish to
understand the tangent space to fibers of the natural morphism
$$
\begin{array}{rccc}
  P : &\Hom(\ell,X) & \to & \Pic(\ell) \\
   & g & \mapsto & g^*(L)
\end{array}
$$
It seems folklore among a handful of experts that the tangent map
$$
TP|_f : \underbrace{T_{\Hom(\ell,X)}|_f}_{= H^0(\ell, f^*(T_X))}
\to \underbrace{T_{\Pic(\ell)}|_{f^*(L)}}_{= H^1(\ell, \O_\ell)}
$$
can be expressed in terms of the first jet sequence of $L$ in the
following way. Dualize Sequence~\eqref{eq:firstJetSeq} and twist by
$L$ to obtain
\begin{equation}  \label{eq:dtJet}
  0 \longrightarrow \O_X \longrightarrow \Jet^1(L)^{\vee}\otimes L
  \longrightarrow T_X \longrightarrow 0.
\end{equation}
The tangent map $TP$ is then the first connecting morphism in the long
exact sequence associated to the $f$-pull-back of \eqref{eq:dtJet},
$$
\cdots \longrightarrow H^0(\ell, f^*(\Jet^1(L)^\vee\otimes L))
\longrightarrow H^0(\ell, f^*(T_X)) \xrightarrow{TP} H^1(\ell,
\O_\ell) \longrightarrow \cdots
$$
For lack of a reference, we will prove the following weaker
statement here which is sufficient for our purposes. More details will
appear in a forthcoming survey.

\begin{thm}\label{thm:ploarized_defo}
  Let $\Delta$ be the unit disc with coordinate $t$ and
  $$
  \begin{array}{rccc}
    f : & \Delta & \to & \Hom(\ell,X) \\
    & t & \mapsto & f_t
  \end{array}
  $$
  be a family of morphisms. Assume that there exists a line bundle
  $H \in \Pic(\ell)$ such that for all $t\in \Delta$ we have $f_t^*(L)
  \cong H$. Then the tangent vector
  $$
  Tf\left(\left.\frac{\partial}{\partial t}\right|_{t=0}\right) \in
  T_{\Hom(\ell,X)}|_{f_0} = H^0(\ell, f_0^*(T_X))
  $$
  is contained in the image of the morphism
  $$
  H^0(\ell, f_0^*(\Jet^1(L)^\vee \otimes L)) \to H^0(\ell,
  f_0^*(T_X))
  $$
  which comes from the dualized and twisted jet
  sequence~\eqref{eq:dtJet}.
\end{thm}

The proof of Theorem~\ref{thm:ploarized_defo} may look rather involved
at first glance, but with the results of Appendix~\ref{app:jets}, its
proof takes little more than a good choice of coordinates on the
projectivized line bundles and an unwinding of the definitions.

\subsubsection*{Proof of Theorem~\ref{thm:ploarized_defo}, Step 1}

Consider the diagram
$$
\xymatrix{
  {\Delta \times \ell} \ar[rrr]^{\mu}_{\text{universal morphism}} 
  \ar[d]_{\text{projection $\pi_2$}}&&& X \\
  {\ell} }
$$
As a first step, we will find convenient coordinates on the
pull-back of the $\P^1$-bundle over $X$,
$$
\mu^* \P(L \oplus\O_X) \cong \P(\mu^*(L) \oplus\O_{\Delta \times
  \ell}).
$$
Since $\Pic(\Delta)=\{e\}$, there exists an isomorphism $\mu^*(L)
\cong \pi_2^*(H)$ which induces an isomorphism
$$
\mu^* \P(L \oplus\O_X) \cong \Delta \times \P(H \oplus \O_\ell).
$$
We use these coordinates to write the base change diagram as
follows:
\begin{equation}\label{diag:jetsandprojs}  
  \xymatrix{ {\Delta \times \P(H \oplus \O_\ell)} \ar[r]^{\tilde \mu}
    \ar[d]_{\tilde \pi}
    & {\P(L \oplus \O_X)} \ar[d]^{\pi}\\
    {\Delta \times \ell} \ar[r]_{\mu} & {X} }
\end{equation}
For convenience of notation, write $Y := \P( L \oplus \O_X)$ and $Y'
:= \Delta \times \P(H \oplus \O_\ell)$.  Let $\Sigma_\ell \subset Y'$
and $\Sigma_X \subset Y$ be the disjoint union of the sections that
come from the direct sum decompositions. It is clear from the
construction that $\tilde \mu(\Sigma_\ell) \subset \Sigma_X$.

\subsubsection*{Proof of Theorem~\ref{thm:ploarized_defo}, Step 2}

Recall from \cite[Rem.~III.9.3.1]{Ha77} that there exists a natural
morphism of sheaves
$$
\alpha : \mu^*\pi_* (T_Y(-\log \Sigma_X)) \to \tilde \pi_* \tilde
\mu^* (T_Y(-\log \Sigma_X))
$$
Although $\mu$ is not flat, we claim the following.

\begin{claim}\label{claim:pushpullpullpush}
  The map $\alpha$ is an isomorphism.
\end{claim}

\begin{proof}
  Because the claim is local on the base, we can assume without loss
  of generality that the locally trivial $\P^1$-bundle $\pi$ is
  actually trivial. For trivial $\P^1$-bundles, however,
  \cite[Prop.~II.3.2(iii)]{Deligne70} shows that the logarithmic
  tangent sheaf decomposes as
  $$
  T_Y(-\log \Sigma_X) \cong \O_Y \oplus \pi^*(T_X).
  $$
  For these two sheaves, Claim~\ref{claim:pushpullpullpush} follows
  easily from the commutativity of Diagram~\eqref{diag:jetsandprojs}
  and from the projection formula.
\end{proof}

\subsubsection*{Proof of Theorem~\ref{thm:ploarized_defo}, End of
  proof}

To avoid confusion, we name the canonical liftings of vector fields on
$\Delta$
\begin{align*}
  \tau_{\rm up} & := \frac{\partial}{\partial t} \in H^0( Y',
  T_{Y'}) \\
  \tau_{\rm down} & := \frac{\partial}{\partial t} \in H^0(\Delta
  \times \ell, T_{\Delta \times \ell})
\end{align*}
In view of Fact~\ref{fact:tgt_vect_to_Hom}, to prove
Theorem~\ref{thm:ploarized_defo}, it suffices to show the following
stronger statement.

\begin{claim}\label{claim:stl}
  The vector field
  $$
  T\mu(\tau_{\rm down}) \in H^0(\Delta \times \ell, \mu^*(T_X))
  $$
  is in the image of 
  $$
  \beta : H^0(\Delta \times \ell, \mu^*(\Jet^1(L)^\vee \otimes L)) \to
  H^0(\Delta \times \ell, \mu^*(T_X)).
  $$
\end{claim}
Use Corollary~\ref{cor:jetsandlogs} and the results of Step (2) to
identify
\begin{align*}
  H^0(\Delta \times \ell, \mu^*(\Jet^1(L)^\vee \otimes L)) & \cong
  H^0( Y', \tilde \mu^* (T_Y(-\log \Sigma_X)))\\
  H^0(\Delta \times \ell, \mu^*(T_X)) & \cong 
  H^0( Y', \tilde \mu^* \pi^*(T_X)).
\end{align*}
These identifications make it easier to write down $\beta$.  Namely,
by Corollary~\ref{cor:jetsandlogs}, $\beta$ becomes nothing but the
pull-back of the tangent map of $\pi$, i.e.~$\beta = \tilde
\mu^*(T\pi)$. Claim~\ref{claim:stl} is thus reformulated as:

\begin{claim}\label{claim:lc}
  The vector field
  $$
  T\tilde \mu(\tau_{\rm up}) \in H^0(Y', \tilde \mu^*(T_X))
  $$
  is in the image of 
  $$
  \tilde \mu^*(T\pi) : H^0(Y', \tilde \mu^*(T_Y(- \log \Sigma_X )))
  \to H^0(Y', \tilde \mu^* \pi^*(T_X)).
  $$
\end{claim}

In this formulation, the proof of Claim~\ref{claim:stl}, and hence of
Theorem~\ref{thm:ploarized_defo}, becomes trivial. The only thing to
note is that
$$
T\tilde \mu (\tau_{\rm up}) \in H^0(Y', \tilde \mu^*(T_Y(- \log
\Sigma_X ))) \subset H^0(Y', \tilde \mu^*(T_Y)).
$$
That, however, follows from the facts that $\tau_{\rm up} \in
H^0(Y', T_{Y'}(-\log \Sigma_\ell))$ and that $\tilde \mu(\Sigma_\ell)
\subset \Sigma_X$. Theorem~\ref{thm:ploarized_defo} is thus shown.
\qed

\subsection{Acknowledgements}

The two appendices would never have appeared without the help of many
people. The author is particularly grateful to Brendan Hassett and
Charles Walter, who suggested to use jets in order to describe the
tangent map $TP$, and to Hubert Flenner and János Kollár for patiently
answering questions by e-mail.

\end{document}